\documentclass{article}
\usepackage[latin1]{inputenc} 
\usepackage{a4wide}           
\usepackage{amssymb}					
\usepackage[arrow, matrix, curve]{xy}
\usepackage{amsmath}					
\usepackage{amsthm}						
\usepackage{bbm}  				  	
\usepackage{paralist}					

\newcommand{\D}{\mathcal D}
\newcommand{\E}{\mathcal E}
\newcommand{\F}{\mathcal F}
\newcommand{\G}{\mathcal G}
\renewcommand{\H}{\mathcal H}
\newcommand{\J}{\mathcal J}
\newcommand{\I}{\mathcal I}
\newcommand{\C}{\mathcal C}
\newcommand{\N}{\mathcal N}
\newcommand{\M}{\mathcal M}
\newtheoremstyle{style}   
  {0.3cm}                  
  {0.3cm}              
  {}                   
  {}                   
  {\normalfont\bfseries}
  {:}{ }      
  {{\normalfont\bfseries \thmname{#1}\thmnumber{ #2}\thmnote{ (#3)}}}
\newtheoremstyle{stylesatz}   
  {0.3cm}                  
  {0.3cm}              
  {\itshape}     
  {}                      
  {\normalfont\bfseries}  
  {:}{ }      
  {{\normalfont\bfseries \thmname{#1}\thmnumber{ #2}\thmnote{ (#3)}}}
\newtheoremstyle{stylebeweis}   
  {0.3cm}              
  {0.3cm}              
  {}                   
  {}                   
  {\normalfont\itshape} 
  {:}{ }      
  {{\normalfont\itshape \thmname{#1}\thmnumber{ #2}\thmnote{ (#3)}}}
\theoremstyle{style}
\newtheorem*{Def}{Definition}
\theoremstyle{stylesatz}
\newtheorem{Satz}{Theorem}\numberwithin{Satz}{section}
\newtheorem{Lemma}[Satz]{Lemma}
\newtheorem{Kor}[Satz]{Corollary}
\newtheorem{Sit}[Satz]{Situation}
\newtheorem*{Satz*}{Theorem}
\theoremstyle{stylebeweis}
\newtheorem*{Bew}{Proof}
\numberwithin{equation}{section}

\begin{document}

\title{Atiyah classes with values\\in the truncated cotangent complex\\\large\bigskip Diploma thesis at the University of Bonn}
\author{Fabian Langholf\\\bigskip\small langholf@gmx.de}
\date{}
\maketitle

\begin{abstract}
We prove an explicit formula for the truncated Atiyah class of a bounded complex of vector bundles. Furthermore, we show that the first truncated Chern class of such a complex only depends on its determinant.
\end{abstract}

\section{Introduction}

Atiyah classes were first introduced in 1957. \textit{Atiyah} associates to vector bundles $\E$ over complex manifolds $X$ elements in $H^1(X,\mathcal{H}om(\E,\E)\otimes\Omega_X)$. They allow him to derive a criterion for the existence of holomorphic connections in \cite[Thm. 2]{Ati}. Furthermore, \textit{Atiyah} remarks that his classes suggest a new way to define the already known Chern classes (\cite[Thm. 6]{Ati}).

Later, it turned out that deformation theory forms an important area for applications of Atiyah classes. If $\E$ is a vector bundle on a scheme $X$, and $X\hookrightarrow Y$ is a closed immersion with nilpotent ideal sheaf, we might be interested in the question whether $\E$ can be extended to a vector bundle on $Y$. \textit{Illusie} gives an extensive answer to questions of this type in his dissertation \cite[Introd. of Chapt. IV, Prop. IV.3.1.8]{Ill}, making essential use of classes that coincide with \textit{Atiyah's} in special cases. For him, the central object is the cotangent complex $\tilde{\mathbb L}_{X|S}$ of a morphism $X\rightarrow S$ of schemes. The outcome shows that this complex is the correct generalization of the cotangent bundle of a smooth morphism. Consequently, his version of the Atiyah class of a vector bundle $\E$ on an $S$-scheme $X$ is an element in ${\rm Ext}^1(\E,\E\otimes\tilde{\mathbb L}_{X|S})$. This leads to the definition of a very general variant of Chern classes, whose properties are studied in \cite[Chapt. V]{Ill}.

\textit{Illusie's} results impress by their overwhelming generality, but they require highly elaborate techniques. Already the definition of his cotangent complex is complicated. 
Fortunately, there is an easier variant of this complex, introduced by \textit{Berthelot} in \cite[Sect. VIII.2]{SGA6}. It is obtained from \textit{Illusie's} cotangent complex by truncation (see \cite[Cor. III.1.2.9.1]{Ill}), thus we call it the truncated cotangent complex ${\mathbb L}_{X|S}$ of a morphism $X\rightarrow S$.

Recently, using the easier complex, \textit{Huybrechts} and \textit{Thomas} managed to prove results similar to \textit{Illusie's} with more elementary methods that also cover the deformation theory of complexes as objects in the derived category. Again, a version of the Atiyah class plays a key role. For a vector bundle $\E$ on an $S$-scheme $X$, this version is an element in ${\rm Ext}^1(\E,\E\otimes{\mathbb L}_{X|S})$. The result motivates us to investigate this truncated Atiyah class in more detail.

In Section 2, we will explain the basic concepts. The third section is devoted to truncated Atiyah classes of complexes $\E$ of vector bundles. Our main result (for the precise formulation see Thm. 3.4) gives an explicit formula via \v Cech resolutions describing these classes:
\begin{Satz*}
The truncated Atiyah class of $\E$ is given by the \v Cech cocycle
\begin{align}
&\left(\left(
M^s_{ik}(\tilde M^s_{kj}\cdot\tilde M^s_{ji}-\tilde M^s_{ki})
\right)_{ijk},\left(
M^s_{ij}\cdot d\tilde M^s_{ji}
\right)_{ij},\right.\\\nonumber&\left.\left(
(-1)^{s+1}M_{ij}^{s+1}(\tilde M_{ji}^{s+1}\cdot\tilde D_i^s-\tilde D_j^s\cdot\tilde M_{ji}^s)
\right)_{ij},\left(
(-1)^{s+1}d\tilde D_i^s
\right)_i,\left(
-\tilde D_i^{s+1}\cdot\tilde D_i^s
\right)_i\right).
\end{align}
\end{Satz*}
Like the other versions of Atiyah classes, the truncated one induces Chern classes. We will discuss the first one in Section 4. It turns out (see Thm. 4.5) the first truncated Chern class of a complex of vector bundles only depends on its determinant:
\begin{Satz*}
For the first truncated Chern class $c_1(\E)\in H^1(X,\mathbb L_X)$, we have
$c_1(\E)=c_1({\rm det}(\E))$.
\end{Satz*}
The result simplifies an argument in \cite{HT} to prove the existence of a perfect obstruction theory for stable pairs as used in \cite{PT}.

\vspace{0.3cm}
\underline{Notation:} We fix a noetherian separated scheme $S$. We are mostly dealing with the category of $S$-schemes, therefore we use the abbreviations $X\times Y:=X\times_S Y$ for fibre products and $\Omega_X:=\Omega_{X|S}$ for relative cotangent sheaves of $S$-schemes $X$ and $Y$.

Rings and algebras are always assumed to be commutative and unitary. If $k$ is a ring, $A$ a $k$-algebra, and $I\subseteq A\otimes_k A$ the ideal of the diagonal, we frequently use the isomorphism $\Omega_{A|k}\cong I/I^2, da\mapsto 1\otimes a-a\otimes 1$ for the module of Kähler differentials --- in the case of morphisms of schemes, we use the analogous isomorphism for the relative cotangent sheaves.

Speaking about sheaves on a scheme $X$, we always mean sheaves of $\mathcal O_X$-modules, whose category is denoted by ${\rm Mod}(X)$. We often consider a sheaf as a complex of sheaves which is concentrated in degree $0$.

Complexes $\E$ of sheaves on a scheme $X$ always have increasing differentials $d^s:\E^s\rightarrow\E^{s+1}$. The complex $\E[1]$ has components $\E[1]^s=\E^{s+1}$ and differentials $d_{\E[1]}^s=-d_\E^{s+1}$. If $\F$ is another complex on $X$, the complex $\E\otimes\F$ with components $(\E\otimes\F)^u=\bigoplus_{s+t=u}\E^s\otimes\F^t$ is induced by the (anti-commuting) bicomplex with differentials $d_\E^s\otimes {\rm id}:\E^s\otimes\F^t\rightarrow\E^{s+1}\otimes\F^{t}$ and $(-1)^s {\rm id}\otimes d_\F^t:\E^s\otimes\F^t\rightarrow\E^{s}\otimes\F^{t+1}$.
The complex $\mathcal{H}om(\E,\F)$ with components $(\mathcal{H}om(\E,\F))^u=\prod_{-s+t=u}\mathcal{H}om(\E^{s},\F^t)$ is induced by the bicomplex with differentials
$\mathcal{H}om(d_\E^{s-1},\F^t):\mathcal{H}om(\E^s,\F^t)\rightarrow\mathcal{H}om(\E^{s-1},\F^t)$ and $(-1)^{-s+t+1}\mathcal{H}om(\E^s,d_\F^{t}):\mathcal{H}om(\E^s,\F^t)\rightarrow\mathcal{H}om(\E^s,\F^{t+1})$.

In the literature, there are different conventions concerning the signs of Atiyah and Chern classes. Basically, there are three choices that have to be made. One of them is the convention concerning the Kähler differentials we introduced above. The other ones appear in the definitions of the classes. Our sign conventions necessarily differ from those of other authors. However, we achieve transparency by avoiding superfluous identifications and by disclosing our choices.

\section{The truncated Atiyah class}

In this section, the truncated cotangent complex $\mathbb L_X$ of a (suitable) $S$-scheme $X$ is to be introduced, and we want to explain how we can associate to each complex $\E$ of sheaves on $X$ a morphism $At(\E):\E\rightarrow\E\otimes^L\mathbb L_X[1]$ in the derived category $\D({\rm Mod}(X))$, called the truncated Atiyah class of $\E$. Furthermore, we will discuss the relationship between the truncated Atiyah class and the classical version $At_{cl}(\E):\E\rightarrow\E\otimes^L\Omega_X[1]$.

\subsection{The truncated cotangent complex}

Our definition of the truncated cotangent complex will not be applicable to all $S$-schemes, but only to those with a smooth ambient space in the sense of the following definition.

\begin{Def}
Let $X$ be an $S$-scheme. A \emph{smooth ambient space} of $X$ is a smooth, separated, and quasi-compact $S$-scheme $U$ with a closed immersion $X\hookrightarrow U$.
\end{Def}

All quasi-projective $S$-schemes have a smooth ambient space because there exists a closed immersion into an open subscheme of a projective space $\mathbb P_S^N$. Furthermore, all $S$-schemes with a smooth ambient space are separated and of finite type over $S$. Hence they are noetherian separated schemes.

\vspace{0.3cm}
Following \textit{Berthelot} (\cite[Sect. VIII.2, Prop. VIII.2.2]{SGA6}), we will define the truncated cotangent complex. The following lemma will make sure that it is well-defined.

\begin{Lemma}[Berthelot]\label{l21}
Let $X$ be an $S$-scheme. Let $\iota_U:X\hookrightarrow U$ and $\iota_V:X\hookrightarrow V$ be smooth ambient spaces, given by ideal sheaves $\J_U$ and $\J_V$. Then the complexes of sheaves $\mathbb L^U_X:=({\J_U/\J_U^2}\stackrel{d}{\rightarrow}\Omega_U|_X )$ and $\mathbb L^V_X:=({\J_V/\J_V^2}\stackrel{d}{\rightarrow}\Omega_V|_X)$ (which are concentrated in degrees \mbox{-1} and 0, and are induced by the conormal sequence) are canonically isomorphic in the derived category $\D({\rm Mod}(X))$.
\end{Lemma}

\begin{Def}
Let $X$ be an $S$-scheme with a smooth ambient space $U$. The \emph{truncated cotangent complex} of $X$ is the complex $\mathbb L_X:=\mathbb L^{U}_X$. It is unique up to canonical isomorphism in the derived category $\D({\rm Mod}(X))$.
\end{Def}

In the situation of Lemma \ref{l21}, the conormal sequences of the two immersions induce morphisms of complexes $\mathbb L^U_X\rightarrow\Omega_X$ and $\mathbb L^V_X\rightarrow\Omega_X$. They are compatible with the canonical isomorphism of the lemma. This follows directly from the definitions. Hence there is a canonical morphism $\mathbb L_X\rightarrow\Omega_X$ in $\D({\rm Mod}(X))$.

If $X$ is smooth over $S$, then this canonical morphism is an isomorphism. We can use $X$ as its own ambient space, and then $\mathbb L_X^X=\Omega_X$.

\subsection{Introduction of the truncated Atiyah class}

In the sequel, we will discuss the truncated Atiyah class, which is introduced by \textit{Huybrechts} and \textit{Thomas} in \cite{HT}. Let $X$ be an $S$-scheme with a smooth ambient space. We want to associate to complexes $\E$ of sheaves on $X$ morphisms $At(\E):\E\rightarrow\E\otimes^L\mathbb L_X[1]$ in a natural way, hence a natural transformation of functors ${\rm id},(\_\otimes^L\mathbb L_X[1]):\D({\rm Mod}(X))\rightrightarrows\D({\rm Mod}(X))$. It is helpful to consider these functors as Fourier--Mukai transforms. From this point of view, the required natural transformation can be defined by giving a single morphism between the Fourier--Mukai kernels $\Delta_{X*}\mathcal O_X$ and $\Delta_{X*}\mathbb L_X[1]$, which we will call universal truncated Atiyah class, and the construction of which is prepared by the following lemma.

\begin{Lemma}[Huybrechts, Thomas]\label{l22}
Let $X$ be a flat $S$-scheme with smooth ambient space $X\hookrightarrow U$, given by the ideal sheaf $\J$. Then the sequence
\begin{equation*}
\begin{xy}
\xymatrix{
0\ar[r]&\Delta_{X*}{\J/\J^2}\ar^\beta[r]&\I_{\Delta_U}|_{X\times X}\ar^\alpha[r]&\mathcal O_{X\times X}\ar^\varepsilon[r]&\Delta_{X*}\mathcal O_X\ar[r]&0
}\end{xy}
\end{equation*}
of sheaves on $X\times X$ is exact. Here $\I_{\Delta_U}$ denotes the ideal of the diagonal on $U\times U$, $\varepsilon$ is the natural surjection, $\alpha$ is the restriction of the inclusion $\I_{\Delta_U}\hookrightarrow\mathcal O_{U\times U}$, and $\beta$ is defined as follows. Let ${\rm Spec}(A)\subseteq U$ be an affine open subset in the preimage of an affine open subset ${\rm Spec}(k)\subseteq S$ with $Y:=X\cap {\rm Spec}(A)={\rm Spec}(A/J)$.
Then $\beta$ maps an element in $\Delta_{X*}{\J/\J^2}(Y\times Y)=J/J^2$, given by a representative $j\in J$, to the class of the element $j\otimes 1-1\otimes j\in I$ in $\I_{\Delta_U}|_{X\times X}(Y\times Y)=I\otimes_{A\otimes_k A}(A/J\otimes_k A/J)$, where $I\subseteq A\otimes_k A$ is the ideal of the diagonal.
\end{Lemma}

\begin{Bew}
The lemma is proved in \cite[Sect. 2.2]{HT}.
\hfill$\Box$
\end{Bew}

The flatness assumption is missing in \cite{HT}, but it is necessary as the following example shows. Consider the situation $S=U={\rm Spec}(\mathbb Z)$, $X={\rm Spec}(\mathbb Z/2\mathbb Z)$. Here, the sequence is of the form $0\rightarrow\mathbb Z/2\mathbb Z\rightarrow 0\rightarrow\mathbb Z/2\mathbb Z\stackrel{{\rm id}}\rightarrow\mathbb Z/2\mathbb Z\rightarrow 0$.

Lemma \ref{l22} is used to define the universal truncated Atiyah class. It makes sure that the morphism between the complexes in the two upper rows of (\ref{univAti1}) is a quasi-isomorphism.

\begin{Def}
Let $X$ be a flat $S$-scheme with smooth ambient space $X\hookrightarrow U$, given by the ideal sheaf $\J$. The \emph{universal truncated Atiyah class} of $X$ is the morphism $At_X:\Delta_{X*}\mathcal O_X\rightarrow\Delta_{X*}\mathbb L_X[1]$ in the derived category $\D({\rm Mod}(X\times X))$ which is given by the diagram
\begin{equation}\label{univAti1}
\begin{xy}
\xymatrix{
&&\Delta_{X*}\mathcal O_X\\
\Delta_{X*}{\J/\J^2}\ar^\beta[r]\ar^{\rm id}[d]&\I_{\Delta_U}|_{X\times X}\ar^\alpha[r]\ar^\pi[d]&\mathcal O_{X\times X}\ar_\varepsilon[u]\\
\Delta_{X*}{\J/\J^2}\ar^/-0.2cm/{-d}[r]&\I_{\Delta_U}/\I_{\Delta_U}^2|_{X\times X}
}\end{xy}
\end{equation}
of complexes of sheaves on $X\times X$ (which are concentrated in degrees $-2$ to $0$).
Here $\pi$ is the natural projection; the other maps have already been introduced in this section (note that $\I_{\Delta_U}/\I_{\Delta_U}^2|_{X\times X}\cong\Delta_{X*}(\Omega_U|_X)$).
\end{Def}

The commutativity of (\ref{univAti1}) follows immediately from the definitions of $\beta$, $\pi$, and $d$. Thus the diagram really describes a morphism in the derived category. It is easy to see that it does not depend on the choice of a smooth ambient space.

We already indicated how the universal class induces the natural transformation we are aiming at. This is made precise by the following definition.

\begin{Def}
Let $X$ be a flat $S$-scheme with a smooth ambient space. Let $p,q:X\times X\rightrightarrows X$ be the natural projections. The \emph{truncated Atiyah class} of $X$ is the natural transformation $At:=Rq_*\circ(\_\otimes^{L}At_X)\circ Lp^*$ between the functors ${\rm id}\cong Rq_*\circ(\_\otimes^{L}\Delta_{X*}\mathcal O_X)\circ Lp^*,(\_\otimes^L\mathbb L_X[1])\cong Rq_*\circ(\_\otimes^{L}\Delta_{X*}\mathbb L_X[1])\circ Lp^*:\D({\rm Mod}(X))\rightrightarrows\D({\rm Mod}(X))$.
\end{Def}

Since this definition is of central importance for us, we repeat it with other words. In the situation of the definition, we consider a complex $\E$ of sheaves on $X$. We are interested in its truncated Atiyah class. Firstly, we have to apply the derived functor $Lp^*:\D({\rm Mod}(X))\rightarrow\D({\rm Mod}(X\times X))$ to $\E$. Here, we can make use of the flatness of $X$, which implies that $Lp^*(\E)=p^*\E$. Next, we apply the derived tensor product with the universal truncated Atiyah class. We get a morphism $p^*\E\otimes^{L}\Delta_{X*}\mathcal O_X\xrightarrow{p^*\E\otimes^{L}At_X} p^*\E\otimes^L\Delta_{X*}\mathbb L_X[1]$. In the general case, it is difficult to control this morphism. However, if $\E$ is a bounded above complex of flat sheaves, it is not necessary to derive the tensor product, and the morphism is given by
\begin{equation*}
\begin{xy}
\xymatrix{
p^*\E\otimes\Delta_{X*}\mathcal O_X&p^*\E\otimes\G\ar[l]\ar[r]&p^*\E\otimes\Delta_{X*}\mathbb L_X[1],
}\end{xy}
\end{equation*}
where $\G$ denotes the complex in the middle row of (\ref{univAti1}), and the left arrow is a quasi-isomorphism.

Finally, the derived functor $Rq_*:\D({\rm Mod}(X\times X))\rightarrow\D({\rm Mod}(X))$ has to be applied. Unlike the last to steps, we can not avoid deriving the functor here. Making use of the general theory of Fourier--Mukai transforms, we can identify $Rq_*(p^*\E\otimes^{L}\Delta_{X*}\mathcal O_X)$ with $\E$ and $Rq_*(p^*\E\otimes^L\Delta_{X*}\mathbb L_X[1])$ with $\E\otimes^L\mathbb L_X[1]$ --- however, a simple and general description of the resulting morphism $At(\E):\E\rightarrow\E\otimes^L\mathbb L_X[1]$ appears inaccessible.

\subsection{Comparison with the classical Atiyah class}

At first sight, the definition of the truncated Atiyah class may appear arbitrary. We will see, however, that it lifts the classical Atiyah class.

\begin{Def}
Let $X$ be a separated $S$-scheme of finite type. The \emph{universal (classical) Atiyah class} of $X$ is the morphism $At_{X,cl}:\Delta_{X*}\mathcal O_X\rightarrow\Delta_{X*}\Omega_X[1]$ in the derived category $\D({\rm Mod}(X\times X))$ which is given by the diagram
\begin{equation}\label{univklasAti}
\begin{xy}
\xymatrix{
&\Delta_{X*}\mathcal O_X\\
\Delta_{X*}\Omega_X\ar[r]\ar^{\rm id}[d]&\mathcal O_{X\times X}/\I_{\Delta_X}^2\ar[u]\\
\Delta_{X*}\Omega_X
}\end{xy}
\end{equation}
of complexes of sheaves on $X\times X$ (which are concentrated in degrees $-1$ and $0$). Here $\I_{\Delta_X}$ denotes the ideal of the diagonal on $X\times X$, and the morphisms without labelling come from the natural short exact sequence (note that $\Delta_{X*}\Omega_X\cong\I_{\Delta_X}/\I_{\Delta_X}^2$).
\end{Def}

In other words, the universal classical Atiyah class is the morphism $\Delta_{X*}\mathcal O_X\rightarrow\Delta_{X*}\Omega_X[1]$ given by the short exact sequence $0\rightarrow\Delta_{X*}\Omega_X\rightarrow\mathcal O_{X\times X}/\I_{\Delta_X}^2\rightarrow\Delta_{X*}\mathcal O_X\rightarrow 0$.

As in the truncated case, the universal class induces the special one.

\begin{Def}
Let $X$ be a separated $S$-scheme of finite type. Let $p,q:X\times X\rightrightarrows X$ be the natural projections. The \emph{(classical) Atiyah class} of $X$ is the natural transformation  $At_{cl}:=Rq_*\circ(\_\otimes^{L}At_{X,cl})\circ Lp^*$ between the functors ${\rm id}\cong Rq_*\circ(\_\otimes^{L}\Delta_{X*}\mathcal O_X)\circ Lp^*,(\_\otimes^L\Omega_X[1])\cong Rq_*\circ(\_\otimes^{L}\Delta_{X*}\Omega_X[1])\circ Lp^*:\D({\rm Mod}(X))\rightrightarrows\D({\rm Mod}(X))$.
\end{Def}

Let $X$ be a flat $S$-scheme with a smooth ambient space, and let $\E$ be a complex of sheaves on $X$. We have defined a classical Atiyah class $At_{cl}(\E):\E\rightarrow\E\otimes^L\Omega_X[1]$ and a truncated one $At(\E):\E\rightarrow\E\otimes^L\mathbb L_X[1]$. The natural map $\mathbb L_X\rightarrow\Omega_X$ induces a morphism $\E\otimes^L\mathbb L_X[1]\rightarrow\E\otimes^L\Omega_X[1]$. The following simple lemma shows that the resulting diagram is commutative.

\begin{Lemma}\label{l23}
Let $X$ be a flat $S$-scheme with a smooth ambient space. Then the natural map $\mathbb L_X\rightarrow\Omega_X$ induces a commutative diagram of natural transformations
\begin{equation*}
\begin{xy}
\xymatrix{
&&(\_\otimes^L\mathbb L_X[1])\ar[dd]\\
{\rm id}\ar^{At}[rru]\ar_{At_{cl}}[rrd]\\
&&(\_\otimes^L\Omega_X[1])\\
}\end{xy}
\end{equation*}
between functors $\D({\rm Mod}(X))\rightarrow\D({\rm Mod}(X))$.
\end{Lemma}

\begin{Bew}
It suffices to prove the commutativity of the diagram
\begin{equation*}
\begin{xy}
\xymatrix{
&&\Delta_{X*}\mathbb L_X[1]\ar[dd]\\
\Delta_{X*}\mathcal O_X\ar^{At_X}[rru]\ar_{At_{X,cl}}[rrd]\\
&&\Delta_{X*}\Omega_X[1].\\
}\end{xy}
\end{equation*}
To this end, we fix a smooth ambient space of $X$. If we denote the complexes in the middle rows of (\ref{univAti1}) and (\ref{univklasAti}) by $\G$ and $\G_{cl}$, then there exists a natural map $\G\rightarrow\G_{cl}$ (in degree $0$ the natural projection, in degree $-1$ the composition $\I_{\Delta_U}|_{X\times X}\stackrel{\pi}\rightarrow\I_{\Delta_U}/\I^2_{\Delta_U}|_{X\times X}\cong\Delta_{X*}(\Omega_U|_X)\rightarrow\Delta_{X*}\Omega_X$) making the diagram
\begin{equation*}
\begin{xy}
\xymatrix{
&\G\ar[r]\ar[dd]\ar[ld]&\Delta_{X*}\mathbb L_X[1]\ar[dd]\\
\Delta_{X*}\mathcal O_X\\
&\G_{cl}\ar[lu]\ar[r]&\Delta_{X*}\Omega_X[1]
}\end{xy}
\end{equation*}
obviously commutative. This implies the assertion of the lemma.\hfill$\Box$
\end{Bew}

\subsection{Alternative description}

In the next section, we will need an alternative description of the universal truncated Atiyah class.

\begin{Lemma}
Let $X$ be a flat $S$-scheme with smooth ambient space $X\hookrightarrow U$, given by the ideal sheaf $\J$. Then, in the commutative diagram
\begin{equation}\label{univAti2}
\begin{xy}
\xymatrix{
&&&\Delta_{X*}\mathcal O_X\ar^{\rm id}[d]\\
\Delta_{X*}{\J/\J^2}\ar^/-0.55cm/{\mbox{\scriptsize$\left(\begin{matrix}-\beta\\{\rm id}\end{matrix}\right)$}}[r]&\I_{\Delta_U}|_{X\times X}\oplus\Delta_{X*}{\J/\J^2}\ar^{\raisebox{0.5cm}{\mbox{\scriptsize$\left(\begin{matrix}-\alpha&0\\\pi&-d\end{matrix}\right)$}}}[r]&\mathcal O_{X\times X}\oplus\I_{\Delta_U}/\I_{\Delta_U}^2|_{X\times X}\ar^/0.7cm/{\mbox{\scriptsize$\left(\begin{matrix}-\varepsilon&0\end{matrix}\right)$}}[r]&\Delta_{X*}\mathcal O_X\\
&\Delta_{X*}{\J/\J^2}\ar_{\mbox{\scriptsize$\left(\begin{matrix}0\\{\rm id}\end{matrix}\right)$}}[u]\ar^{-d}[r]&\I_{\Delta_U}/\I_{\Delta_U}^2|_{X\times X},\ar_{\mbox{\scriptsize$\left(\begin{matrix}0\\{\rm id}\end{matrix}\right)$}}[u]
}\end{xy}
\end{equation}
the morphism between the complexes in the two lower rows (which are concentrated in degrees $-3$ to $0$) is a quasi-isomorphism. The morphism $\Delta_{X*}\mathcal O_X\rightarrow\Delta_{X*}\mathbb L_X[1]$ given by the diagram is the universal truncated Atiyah class of $X$.
\end{Lemma}

\begin{Bew}
Once more, we denote the complex in the middle row of (\ref{univAti1}) by $\G$. The complex $\G_0$ in the middle row of (\ref{univAti2}) is the mapping cone of the morphism $\mbox{\scriptsize$\left(\begin{matrix}-\varepsilon\\({\rm id},\pi)\end{matrix}\right)$}:\G\rightarrow\Delta_{X*}\mathcal O_X\oplus\Delta_{X*}\mathbb L_X[1]$. Since $\G\stackrel{\varepsilon}\rightarrow\Delta_{X*}\mathcal O_X$ is a quasi-isomorphism, this holds also true for the induced map $\Delta_{X*}\mathbb L_X[1]\rightarrow\G_0$.
By construction, the composition $\G\rightarrow\Delta_{X*}\mathcal O_X\oplus\Delta_{X*}\mathbb L_X[1]\rightarrow\G_0$ is the zero map, hence the diagram
\begin{equation*}
\begin{xy}
\xymatrix{
\G\ar^\varepsilon[r]\ar[d]_{({\rm id},\pi)}&\Delta_{X*}\mathcal O_X\ar[d]\\
\Delta_{X*}\mathbb L_X[1]\ar[r]&\G_0,
}\end{xy}
\end{equation*}
is commutative, from which the assertion follows.
\hfill$\Box$
\end{Bew}

\section{The truncated Atiyah class of a complex of vector bundles}

In the previous section, we introduced the truncated Atiyah class $At(\E):\E\rightarrow\E\otimes^L\mathbb L_X[1]$ of a complex $\E$ of sheaves on a flat $S$-scheme $X$ with a smooth ambient space. Our definition via a universal class has the advantage that it can be applied to all complexes of sheaves, and that its naturality is obvious. However, it is not suitable for explicit calculations. Without changing our point of view, we will even fail to prove elementary results about the first truncated Chern class, which will be discussed in the next section as an application of the truncated Atiyah class.

Thus we will have to give an alternative, more concrete description of the truncated Atiyah class. We will make use of \v Cech resolutions. For the increase of concreteness, we have to accept a loss of generality. The methods of this section only allow us to handle bounded complexes of vector bundles.

\subsection{Preliminaries}

First, we fix the setting of this section.

\begin{Sit}\label{sit}
Let $X$ be a flat $S$-scheme with smooth ambient space $\iota:X\hookrightarrow U$ (in particular, $U$ is separated), given by the ideal sheaf $\J$. We denote the natural projections by $p,q:X\times X\rightrightarrows X$.

Let $\E$ be a bounded complex of locally free coherent sheaves on $X$. Let $U=\bigcup_\Gamma U_i$ be a finite cover of $U$ by affine open subsets $U_i$ such that the components $\E^s$ of $\E$ are free on all $X_i:=U_i\cap X$, and such that $U_i$ maps into an affine open subset $S_i$ of $S$. We assume that the index set $\Gamma$ is strictly ordered.

For non-empty subsets $\Lambda $ von $\Gamma$, we introduce the notation $U_\Lambda :=\bigcap_\Lambda  U_i$, $X_\Lambda :=\bigcap_\Lambda  X_i$, and $S_\Lambda :=\bigcap_\Lambda  S_i$.
Since $U$, $X$, and $S$ are separated, $U_\Lambda ={\rm Spec}(A_\Lambda )$, $X_\Lambda ={\rm Spec}(A_\Lambda /J_\Lambda )$, and $S_\Lambda ={\rm Spec}(k_\Lambda )$ are affine.

Let $\varphi_i^s:\E^s|_{X_i}\stackrel{\sim}{\rightarrow}\mathcal O_{X_i}^{m_i^s}$ be trivializations and $M^s_{ij}:=\varphi_i^s\circ(\varphi_j^s)^{-1}|_{X_{ij}}\in {\rm GL}(m_i^s,A_{ij}/J_{ij})$ the corresponding transition maps. (Note that in the case $X_{ij}\neq\emptyset$, we always have $m_i^s=m_j^s$.) For each ordered pair $(i,j)$, we choose a lift $\tilde M^s_{ij}\in {\rm Mat}(m^s_i\times m_i^s,A_{ij})$ von $M^s_{ij}$. Similarly, we define $D^s_i:=\varphi_i^{s+1}\circ d_\E^s\circ(\varphi_i^s)^{-1}\in {\rm Mat}(m_i^{s+1}\times m_i^s,A_i/J_i)$ and choose a lift $\tilde D^s_i\in {\rm Mat}(m_i^{s+1}\times m_i^s,A_i)$.
\end{Sit}

Now, we want to bring \v Cech resolutions into play. Let $Y$ be a noetherian separated scheme and $Y=\bigcup Y_i$ a finite affine open cover of $Y$ with strictly ordered index set. For a quasi-coherent sheaf $\F$ on $Y$, we denote its \v Cech resolution (with respect to the chosen cover, see \cite[Lem. III.4.2]{Ha}) by $\C(\F)$. Thus for natural numbers $r$, we have $\C^r(\F)=\bigoplus_\Lambda i_{\Lambda *}(\F|_{\bigcap_\Lambda Y_i})$, where the sum runs over strictly increasing sequences $\Lambda$ of length $r+1$, and $i_\Lambda:\bigcap_\Lambda Y_i\hookrightarrow Y$ denotes the natural inclusion.

If $\H$ is a bounded complex of quasi-coherent sheaves on $Y$, we write $\C(\H)$ for the total complex associated with the (anti-commuting) bicomplex $\C^r(\H^s)$ with differentials $\check{d}:\C^r(\H^s)\rightarrow\C^{r+1}(\H^s)$ and $(-1)^rd_\H:\C^r(\H^s)\rightarrow\C^r(\H^{s+1})$. If we consider $\H$ as a bicomplex concentrated in column $0$, the natural maps $\H^s\rightarrow\C^0(\H^s)$ induce a morphism of bicomplexes and thus a morphism between their total complexes $\H\rightarrow\C(\H)$. For the following, see e.g. \cite[Thm. 12.5.4]{KS}.

\begin{Lemma}\label{l32}
Let $\H$ be a bounded complex of quasi-coherent sheaves on a noetherian separated scheme $Y$. Let $Y=\bigcup Y_i$ be a finite affine open cover of $Y$ with strictly ordered index set. Then the natural morphism $\H\rightarrow\C(\H)$ is a quasi-isomorphism.
\end{Lemma}

Next, we will define a morphism $\E\rightarrow\E\otimes\mathbb L_X[1]$ in the derived category. It will turn out that it coincides with the truncated Atiyah class of $\E$. It follows from Lemma \ref{l32} that $\E\otimes\mathbb L_X[1]$ is quasi-isomorphic to $\E\otimes\C(\mathbb L_X[1])$ in a natural way, where the \v Cech resolution is formed with respect to the cover of $X$ we fixed above. Thus, for the construction of our morphism, it suffices to give a morphism of complexes $\E\rightarrow\E\otimes\C(\mathbb L_X[1])$.

The sheaf $(\E\otimes\C(\mathbb L_X[1]))^s$ has the summand $\E^{s}\otimes\C^2(\J/\J^2)\oplus\E^{s}\otimes\C^1(\Omega_U|_X)\oplus\E^{s+1}\otimes\C^1(\J/\J^2)\oplus\E^{s+1}\otimes\C^0(\Omega_U|_X)\oplus\E^{s+2}\otimes\C^0(\J/\J^2)$. If $\F$ is a sheaf on $X$, then (with the notation introduced above) a morphism  $\E^s\rightarrow\E^t\otimes\C^r(\F)=\E^t\otimes\bigoplus_\Lambda i_{\Lambda *}(\F|_{X_\Lambda})\cong\bigoplus_\Lambda i_{\Lambda *}(\E^t|_{X_\Lambda}\otimes\F|_{X_\Lambda})$ is, by adjunction, given by maps $\E^s|_{X_{\Lambda}}\rightarrow\E^t|_{X_\Lambda}\otimes\F|_{X_{\Lambda}}$. If we use the restrictions of $\varphi^s_{min(\Lambda)}$ and $\varphi^t_{min(\Lambda)}$ as trivializations on $X_\Lambda$ (what we will always do in the sequel), then these maps can be considered as matrices with entries in $\F(X_\Lambda)$. Hence, for the definition of a morphism $\E^s\rightarrow(\E\otimes\C(\mathbb L_X[1]))^s$ (factorizing over the summand mentioned above), it suffices to give matrices in ${\rm Mat}(m_i^s\times m_i^s,J_{ijk}/J_{ijk}^2)$, ${\rm Mat}(m_i^s\times m_i^s,\Omega_{A_{ij}|k_{ij}}\otimes_{A_{ij}}A_{ij}/J_{ij})$, ${\rm Mat}(m_i^{s+1}\times m_i^s,J_{ij}/J_{ij}^2)$, ${\rm Mat}(m_i^{s+1}\times m_i^s,\Omega_{A_{i}|k_{i}}\otimes_{A_{i}}A_{i}/J_{i})$, and ${\rm Mat}(m_i^{s+2}\times m_i^s,J_{i}/J_{i}^2)$.

With these considerations in mind, we define a morphism $\E^s\rightarrow(\E\otimes\C(\mathbb L_X[1]))^s$ by
\begin{align}\label{FormelAtiyah}
&\left(\left(
M^s_{ik}(\tilde M^s_{kj}\cdot\tilde M^s_{ji}-\tilde M^s_{ki})
\right)_{ijk},\left(
M^s_{ij}\cdot d\tilde M^s_{ji}
\right)_{ij},\right.\\\nonumber&\left.\left(
(-1)^{s+1}M_{ij}^{s+1}(\tilde M_{ji}^{s+1}\cdot\tilde D_i^s-\tilde D_j^s\cdot\tilde M_{ji}^s)
\right)_{ij},\left(
(-1)^{s+1}d\tilde D_i^s
\right)_i,\left(
-\tilde D_i^{s+1}\cdot\tilde D_i^s
\right)_i\right).
\end{align}
Here and in the sequel, we drop restrictions to open subsets from the notation. We denote the universal derivation by $d:A_{\Lambda}\rightarrow\Omega_{A_{\Lambda}|k_{\Lambda}}$, which has to be applied separately to the entries of the matrices. The equations $M^s_{kj}\cdot M^s_{ji}=M^s_{ki}$, $M_{ji}^{s+1}\cdot D_i^s=D_j^s\cdot M_{ji}^s$, and $D_i^{s+1}\cdot D_i^s=0$ imply that $\tilde M^s_{kj}\cdot\tilde M^s_{ji}-\tilde M^s_{ki}$, $\tilde M_{ji}^{s+1}\cdot\tilde D_i^s-\tilde D_j^s\cdot\tilde M_{ji}^s$, and $\tilde D_i^{s+1}\cdot\tilde D_i^s$ really are matrices with entries in $J_{\Lambda}$.

We want to motivate the formula (\ref{FormelAtiyah}) in the case of a vector bundle $\E$ with transition maps $M_{ij}$. In this situation, only the first two of the five parts of the formula play a role. It is known that the classical Atiyah class of $\E$ is given by the \v Cech cocycle \mbox{$(M_{ij}\cdot d M_{ji})_{ij}:\E\rightarrow\E\otimes\C^1(\Omega_X)$} --- we will discuss this formula (\ref{FormelklasAti}) later. In view of Lemma \ref{l23}, it appears sensible to search for a morphism $\E\rightarrow\E\otimes\C^1(\mathbb L_X)$ lifting this cocyle and to hope that it describes the truncated class. If $\E$ extends to a vector bundle on $U$, then the transition maps $\tilde M_{ij}$ of the extension yield natural lifts of those of $\E$ (if we choose trivializations of $\E$ which are restrictions of trivializations of the extension). The resulting map $(M_{ij}\cdot d \tilde M_{ji})_{ij}:\E\rightarrow\E\otimes\C^1(\Omega_U|_X)$ really gives a morphism of complexes $\E\rightarrow\E\otimes\C(\mathbb L_X[1])$ solving our problem. If $\E$, however, cannot be extended, then the choice of lifts $\tilde M_{ij}$ cannot be carried out in a natural way. The term $(M_{ik}(\tilde M_{kj}\cdot\tilde M_{ji}-\tilde M_{ki}))_{ijk}$ in (\ref{FormelAtiyah}) is needed as compensation --- it makes sure that a morphism of complexes $\E\rightarrow\E\otimes\C(\mathbb L_X[1])$ arises, as our next lemma will show.

After this interlude, we return to the more general situation \ref{sit} of a complex $\E$.

\begin{Lemma}\label{l33}
In situation \ref{sit}, the morphisms (\ref{FormelAtiyah}) define a morphism of complexes $\E\rightarrow\E\otimes\C(\mathbb L_X[1])$ and thus a morphism $\E\rightarrow\E\otimes\mathbb L_X[1]$ in the derived category $\D({\rm Mod}(X))$.
\end{Lemma}

\begin{Bew}
We have to make sure that the diagrams
\begin{equation*}
\begin{xy}
\xymatrix{
\E^s\ar[r]\ar[ddd]^{d_\E}&{\begin{matrix}\E^{s}\otimes\C^2(\J/\J^2)\oplus\E^{s}\otimes\C^1(\Omega_U|_X)\oplus\E^{s+1}\otimes\C^1(\J/\J^2)\oplus\\\E^{s+1}\otimes\C^0(\Omega_U|_X)\oplus\E^{s+2}\otimes\C^0(\J/\J^2)\end{matrix}}
\ar[ddd]^/-0.2cm/{\mbox{\scriptsize$\begin{pmatrix}\xi \check d&0&0&0&0\\-\xi d&\xi\check d&0&0&0\\d_\E&0&-\xi\check d&0&0\\0&d_\E&-\xi d&-\xi\check d&0\\0&0&d_\E&0&\xi\check d\\0&0&0&d_\E&-\xi d\\0&0&0&0&d_\E\end{pmatrix}$}}
\\\\\\
\E^{s+1}\ar[r]&{\begin{matrix}\E^{s}\otimes\C^3(\J/\J^2)\oplus\E^{s}\otimes\C^2(\Omega_U|_X)\oplus\E^{s+1}\otimes\C^2(\J/\J^2)\oplus\\\E^{s+1}\otimes\C^1(\Omega_U|_X)\oplus\E^{s+2}\otimes\C^1(\J/\J^2)\oplus\\\E^{s+2}\otimes\C^0(\Omega_U|_X)\oplus\E^{s+3}\otimes\C^0(\J/\J^2)
\end{matrix}}
}\end{xy}
\end{equation*}
commute, where the rows are given by (\ref{FormelAtiyah}). Here, $d$ is induced by the universal derivation, and we use the abbreviation $\xi:=(-1)^s$.

Corresponding to the decomposition of the sheaf in the right lower corner of the diagram into seven summands, we have to compare seven pairs of morphisms $\E^s\rightarrow\E^t\otimes\C^r(\F)$ (with $\F=\J/\J^2$ or $\F=\Omega_U|_X$). To this end, we fix strictly increasing sequences of indices $\Lambda=(i,j,k,l)$, $\Lambda=(i,j,k)$, $\Lambda=(i,j)$ or $\Lambda=(i)$ and prove the equality of the two maps $\E^s|_{X_{\Lambda}}\rightarrow\E^t|_{X_\Lambda}\otimes\F|_{X_{\Lambda}}$ by a calculation with the corresponding matrices.

For the first summand, we get by definition of the \v Cech differential
\begin{align*}
&\xi\check d((M^s_{ik}(\tilde M^s_{kj}\cdot\tilde M^s_{ji}-\tilde M^s_{ki}))_{ijk})_{ijkl}=\\
&\xi[M^s_{ij}(M^s_{jl}(\tilde M^s_{lk}\cdot\tilde M^s_{kj}-\tilde M^s_{lj}))M^s_{ji}-M^s_{il}(\tilde M^s_{lk}\cdot\tilde M^s_{ki}-\tilde M^s_{li})+\\&M^s_{il}(\tilde M^s_{lj}\cdot\tilde M^s_{ji}-\tilde M^s_{li})-M^s_{ik}(\tilde M^s_{kj}\cdot\tilde M^s_{ji}-\tilde M^s_{ki})]=\\&
\xi[\tilde M^s_{il}\cdot\tilde M^s_{lk}\cdot\tilde M^s_{kj}\cdot\tilde M^s_{ji}-\tilde M^s_{il}\cdot\tilde M^s_{lj}\cdot\tilde M^s_{ji}-\tilde M^s_{il}\cdot\tilde M^s_{lk}\cdot\tilde M^s_{ki}+\tilde M^s_{il}\cdot\tilde M^s_{li}+\\&\tilde M^s_{il}\cdot\tilde M^s_{lj}\cdot\tilde M^s_{ji}-\tilde M^s_{il}\cdot\tilde M^s_{li}-\tilde M^s_{ik}\cdot\tilde M^s_{kj}\cdot\tilde M^s_{ji}+\tilde M^s_{ik}\cdot\tilde M^s_{ki}]=\\&
\xi(\tilde M^s_{il}\cdot\tilde M^s_{lk}-\tilde M^s_{ik})\cdot(\tilde M^s_{kj}\cdot\tilde M^s_{ji}-\tilde M^s_{ki})=0.
\end{align*}
The last equality results from the fact that the product of two matrices with entries in $J_{ijkl}$ has entries in $J_{ijkl}^2$.

For the second summand, we compute
\begin{align*}
&-\xi d(M^s_{ik}(\tilde M^s_{kj}\cdot\tilde M^s_{ji}-\tilde M^s_{ki}))+\xi\check d((M^s_{ij}\cdot d\tilde M^s_{ji})_{ij})_{ijk}=\\
&\xi[-M^s_{ik}\cdot d(\tilde M^s_{kj}\cdot\tilde M^s_{ji}-\tilde M^s_{ki})+M^s_{ij}(M^s_{jk}\cdot d\tilde M^s_{kj})M^s_{ji}-M^s_{ik}\cdot d\tilde M^s_{ki}+M^s_{ij}\cdot d\tilde M^s_{ji}]=\\&
\xi[-M^s_{ik}\cdot d\tilde M^s_{kj}\cdot M^s_{ji}-M^s_{ik}\cdot M^s_{kj}\cdot d\tilde M^s_{ji}+M^s_{ik}\cdot d\tilde M^s_{ki}+\\&M^s_{ik}\cdot d\tilde M^s_{kj}\cdot M^s_{ji}-M^s_{ik}\cdot d\tilde M^s_{ki}+M^s_{ij}\cdot d\tilde M^s_{ji}]=0.
\end{align*}

For the third summand, we get (using $D^s_i\cdot M^s_{ik}=M^{s+1}_{ik}\cdot D^s_k$ for the second equality)
\begin{align*}
&d_\E(M^s_{ik}(\tilde M^s_{kj}\cdot\tilde M^s_{ji}-\tilde M^s_{ki}))-\xi\check d((-\xi M_{ij}^{s+1}(\tilde M_{ji}^{s+1}\cdot\tilde D_i^s-\tilde D_j^s\cdot\tilde M_{ji}^s))_{ij})_{ijk}=\\&
D^s_i\cdot M^s_{ik}(\tilde M^s_{kj}\cdot\tilde M^s_{ji}-\tilde M^s_{ki})+M_{ij}^{s+1}(M_{jk}^{s+1}(\tilde M_{kj}^{s+1}\cdot\tilde D_j^s-\tilde D_k^s\cdot\tilde M_{kj}^s))M_{ji}^s-\\&M_{ik}^{s+1}(\tilde M_{ki}^{s+1}\cdot\tilde D_i^s-\tilde D_k^s\cdot\tilde M_{ki}^s)+M_{ij}^{s+1}(\tilde M_{ji}^{s+1}\cdot\tilde D_i^s-\tilde D_j^s\cdot\tilde M_{ji}^s)=\\&
\tilde M^{s+1}_{ik}\cdot\tilde D^s_k\cdot\tilde M^s_{kj}\cdot\tilde M^s_{ji}-\tilde M^{s+1}_{ik}\cdot\tilde D^s_k\cdot\tilde M^s_{ki}+\tilde M_{ik}^{s+1}\cdot\tilde M_{kj}^{s+1}\cdot\tilde D_j^s\cdot\tilde M_{ji}^s-\tilde M_{ik}^{s+1}\cdot\tilde D_k^s\cdot\tilde M_{kj}^s\cdot\tilde M_{ji}^s-\\&\tilde M_{ik}^{s+1}\cdot\tilde M_{ki}^{s+1}\cdot\tilde D_i^s+\tilde M_{ik}^{s+1}\cdot\tilde D_k^s\cdot\tilde M_{ki}^s+\tilde M_{ik}^{s+1}\cdot\tilde M_{kj}^{s+1}\cdot\tilde M_{ji}^{s+1}\cdot\tilde D_i^s-\tilde M_{ik}^{s+1}\cdot\tilde M_{kj}^{s+1}\cdot\tilde D_j^s\cdot\tilde M_{ji}^s=\\&
M_{ik}^{s+1}(\tilde M_{kj}^{s+1}\cdot\tilde M_{ji}^{s+1}-\tilde M_{ki}^{s+1})D_i^s=
(M^{s+1}_{ik}(\tilde M^{s+1}_{kj}\cdot\tilde M^{s+1}_{ji}-\tilde M^{s+1}_{ki}))d_\E.
\end{align*}

For the fourth summand, we get (using $D_i^s\cdot M^s_{ij}=M_{ij}^{s+1}\cdot D_j^s$ for the last but one equality)
\begin{align*}
&d_\E(M^s_{ij}\cdot d\tilde M^s_{ji})-\xi d(-\xi M_{ij}^{s+1}(\tilde M_{ji}^{s+1}\cdot\tilde D_i^s-\tilde D_j^s\cdot\tilde M_{ji}^s))-\xi\check d((-\xi d\tilde D_i^s)_i)_{ij}=\\&
D_i^s\cdot M^s_{ij}\cdot d\tilde M^s_{ji}+M_{ij}^{s+1}\cdot d(\tilde M_{ji}^{s+1}\cdot\tilde D_i^s-\tilde D_j^s\cdot\tilde M_{ji}^s)+M_{ij}^{s+1}\cdot d\tilde D_j^s\cdot M_{ji}^s-d\tilde D_i^s=\\&
D_i^s\cdot M^s_{ij}\cdot d\tilde M^s_{ji}+M_{ij}^{s+1}\cdot d\tilde M_{ji}^{s+1}\cdot D_i^s+M_{ij}^{s+1}\cdot M_{ji}^{s+1}\cdot d\tilde D_i^s-\\&M_{ij}^{s+1}\cdot d\tilde D_j^s\cdot M_{ji}^s-M_{ij}^{s+1}\cdot D_j^s\cdot d\tilde M_{ji}^s+M_{ij}^{s+1}\cdot d\tilde D_j^s\cdot M_{ji}^s-d\tilde D_i^s=\\&
M_{ij}^{s+1}\cdot d\tilde M_{ji}^{s+1}\cdot D_i^s=(M^{s+1}_{ij}\cdot d\tilde M^{s+1}_{ji})d_\E.
\end{align*}

For the fifth summand, we get
\begin{align*}
&d_\E(-\xi M_{ij}^{s+1}(\tilde M_{ji}^{s+1}\cdot\tilde D_i^s-\tilde D_j^s\cdot\tilde M_{ji}^s))+\xi\check d((-\tilde D_i^{s+1}\cdot\tilde D_i^s)_i)_{ij}=\\&
\xi[-D_i^{s+1}\cdot M_{ij}^{s+1}(\tilde M_{ji}^{s+1}\cdot\tilde D_i^s-\tilde D_j^s\cdot\tilde M_{ji}^s)-M_{ij}^{s+2}\cdot\tilde D_j^{s+1}\cdot\tilde D_j^s\cdot M_{ji}^s+\tilde D_i^{s+1}\cdot\tilde D_i^s]=\\&
\xi[-M_{ij}^{s+2}\cdot D_j^{s+1}(\tilde M_{ji}^{s+1}\cdot\tilde D_i^s-\tilde D_j^s\cdot\tilde M_{ji}^s)-M_{ij}^{s+2}\cdot\tilde D_j^{s+1}\cdot\tilde D_j^s\cdot M_{ji}^s+\tilde D_i^{s+1}\cdot\tilde D_i^s]=\\&
\xi[-\tilde M_{ij}^{s+2}\cdot\tilde D_j^{s+1}\cdot\tilde M_{ji}^{s+1}\cdot\tilde D_i^s+M_{ij}^{s+2}\cdot M_{ji}^{s+2}\cdot\tilde D_i^{s+1}\cdot\tilde D_i^s]=\\&
\xi M_{ij}^{s+2}(\tilde M_{ji}^{s+2}\cdot\tilde D_i^{s+1}-\tilde D_j^{s+1}\cdot\tilde M_{ji}^{s+1})D_i^s=(\xi M_{ij}^{s+2}(\tilde M_{ji}^{s+2}\cdot\tilde D_i^{s+1}-\tilde D_j^{s+1}\cdot\tilde M_{ji}^{s+1}))d_\E.
\end{align*}

For the sixth summand, we get
\begin{align*}
&d_\E(-\xi d\tilde D_i^s)-\xi d(-\tilde D_i^{s+1}\cdot\tilde D_i^s)=\xi[-D_i^{s+1}\cdot d\tilde D_i^s+d\tilde D_i^{s+1}\cdot D_i^s+D_i^{s+1}\cdot d\tilde D_i^s)]=\\&
\xi d\tilde D_i^{s+1}\cdot D_i^s=(\xi d\tilde D_i^s)d_\E.
\end{align*}

Finally, for the seventh summand, we have
\begin{align*}
&d_\E(-\tilde D_i^{s+1}\cdot\tilde D_i^s)=-D_i^{s+2}\cdot\tilde D_i^{s+1}\cdot\tilde D_i^s=-\tilde D_i^{s+2}\cdot\tilde D_i^{s+1}\cdot D_i^s=(-\tilde D_i^{s+2}\cdot\tilde D_i^{s+1})d_\E.
\end{align*}
\hfill$\Box$
\end{Bew}

\subsection{Concrete description}

We have constructed two morphisms $\E\rightrightarrows\E\otimes\mathbb L_X[1]$, the rather abstract truncated Atiyah class and the concrete morphism of Lemma \ref{l33}. Our main result shows that they coincide.

\begin{Satz}
In situation \ref{sit}, the morphism $\E\rightarrow\E\otimes\mathbb L_X[1]$ in the derived category $\D({\rm Mod}(X))$ given by the formula (\ref{FormelAtiyah}) is the truncated Atiyah class of $\E$.
\end{Satz}

\begin{Bew}
Let $(X\times X)\backslash\Delta_X=\bigcup_{\Gamma'}V_i$ be a finite cover by affine open sets. For $i\in\Gamma$, we define $V_i:=X_i\times X_i\subseteq X\times X$. Then $X\times X=\bigcup_{\Gamma\sqcup\Gamma'}V_i$ is a finite cover by affine open sets (note that $X_i\times X_i={\rm Spec}(A_i/J_i\otimes_{k_i}A_i/J_i)$). We extend the strict order on $\Gamma$ to a strict order on $\Gamma\sqcup\Gamma'$. For non-empty subsets $\Lambda$ of $\Gamma\sqcup\Gamma'$, we define $V_\Lambda:=\bigcap_\Lambda V_i$, and we denote the natural inclusion by $i_\Lambda:V_\Lambda\rightarrow X\times X$.

We write $\G_0$ for the complex in the middle row of (\ref{univAti2}), and we consider its \v Cech resolution $\C(\G_0)$ with respect to the chosen cover of $X\times X$. We have $\C^r(\G_0^s)=\bigoplus_\Lambda i_{\Lambda *}(\G_0^s|_{V_\Lambda})$, where the sum runs over strictly increasing sequences $\Lambda$ of length $r+1$. We form a subsheaf $\G_1^{rs}$ of $\C^r(\G_0^s)$ by considering only those summands with $\Lambda\cap\Gamma'\neq\emptyset$. The totality of these subsheaves defines a subcomplex of $\C(\G_0)$, denoted by $\G_1$. Further, let $\G:=\C(\G_0)/\G_1$ be the corresponding quotient.

From the condition $\Lambda\cap\Gamma'\neq\emptyset$, it follows that $V_\Lambda\cap\Delta_X=\emptyset$. Thus sheaves which appear in the complex $\G_0$, but are concentrated on the diagonal, do not contribute to $\G_1$. More explicit, $\G_1$ is only composed of the summands $G_1^{r,-1}=\bigoplus_{\Lambda\cap\Gamma'\neq\emptyset}i_{\Lambda *}(\mathcal O_{X\times X}|_{V_\Lambda})$ and $G_1^{r,-2}=\bigoplus_{\Lambda\cap\Gamma'\neq\emptyset}i_{\Lambda *}(\I_{\Delta_U}|_{V_\Lambda})$.
Thus if we write $e:=|\Gamma\sqcup\Gamma'|$, $\G_1$ is given by
\begin{align*}
&\G_1^{0,-2}\xrightarrow{\mbox{\scriptsize$\left(\begin{matrix}\check d\\-\alpha\end{matrix}\right)$}}\G_1^{1,-2}\oplus\G_1^{0,-1}\xrightarrow{\mbox{\scriptsize$\left(\begin{matrix}\check d&0\\\alpha&\check d\end{matrix}\right)$}}\G_1^{2,-2}\oplus\G_1^{1,-1}\rightarrow\cdots\rightarrow
\G_1^{r+1,-2}\oplus\G_1^{r,-1}\\&
\xrightarrow{\mbox{\scriptsize$\left(\begin{matrix}\check d&0\\(-1)^r\alpha&\check d\end{matrix}\right)$}}\G_1^{r+2,-2}\oplus\G_1^{r+1,-1}\rightarrow\cdots\rightarrow
\G_1^{e-1,-2}\oplus\G_1^{e-2,-1}\xrightarrow{\mbox{\scriptsize$\left(\begin{matrix}(-1)^{e-2}\alpha&\check d\end{matrix}\right)$}}\G_1^{e-1,-1}.
\end{align*}
Apart from the diagonal, however, $\alpha$ is an isomorphism, and hence it induces isomorphisms $\G_1^{r,-2}\rightarrow\G_1^{r,-1}$. It follows that $\G_1$ is an exact complex. Thus $\C(\G_0)\rightarrow\G$ is a quasi-isomorphism.

We have constructed the commutative diagram
\begin{equation*}
\begin{xy}
\xymatrix{
\Delta_{X*}\mathcal O_X\ar[r]&\G_0\ar[d]&\Delta_{X*}\mathbb L_X[1]\ar[l]\ar[d]\\
&\C(\G_0)\ar[d]&\C(\Delta_{X*}\mathbb L_X[1])\ar[l]\\
&\G
}\end{xy}
\end{equation*}
of complexes of sheaves on $X\times X$. Here, all morphisms except the left one are quasi-isomorphisms, and the first row describes the truncated Atiyah class of $X$. If we apply the functor $q_*(p^*\E\otimes\_)$, we obtain the diagram
\begin{equation}\label{vekt0}
\begin{xy}
\xymatrix{
\E\ar[r]&q_*(p^*\E\otimes\G_0)\ar[d]&\E\otimes\mathbb L_X[1]\ar[d]\ar[l]\\
&q_*(p^*\E\otimes\C(\G_0))\ar[d]&q_*(p^*\E\otimes\C(\Delta_{X*}\mathbb L_X[1]))\ar[l]\\
&q_*(p^*\E\otimes\G)\ar[d]\\
&Rq_*(p^*\E\otimes\G).
}\end{xy}
\end{equation}
of complexes of sheaves on $X$.
The map $\E\otimes\mathbb L_X[1]\rightarrow Rq_*(p^*\E\otimes\G)$ is a quasi-isomorphism, and the truncated Atiyah class of $\E$ is given by $\E\rightarrow Rq_*(p^*\E\otimes\G)$.

Since $\Delta_{X*}\mathbb L_X[1]$ is concentrated on the diagonal, we have $q_*(p^*\E\otimes\C(\Delta_{X*}\mathbb L_X[1]))=q_*(p^*\E\otimes\Delta_{X*}\C(\mathbb L_X[1]))=\E\otimes\C(\mathbb L_X[1])$ --- note that we compare \v Cech complexes on $X\times X$ (with the cover indexed by $\Gamma\sqcup\Gamma'$) with \v Cech complexes on $X$ (with the cover indexed by $\Gamma$). Thus the morphism $\E\rightarrow\E\otimes\C(\mathbb L_X[1])$ defined in (\ref{FormelAtiyah}) can be inserted into our diagram, and it induces a second morphism $\E\rightarrow q_*(p^*\E\otimes\G)$ (besides the one given by (\ref{vekt0})). We will show that these two morphisms are homotopic. Then it follows that the induced morphisms $\E\rightarrow Rq_*(p^*\E\otimes\G)$ are homotopic as well, and this implies the assertion of the theorem.

Throughout the rest of the proof, we will frequently consider sheaves on $X\times X$ that are pushed forward from $X$ via the diagonal map. In this case, we will omit the term $\Delta_{X*}$ in order to simplify the notation.

A morphism $\E\rightarrow q_*(p^*\E\otimes\G)$ is given by maps $\E^s\rightarrow q_*(p^*\E\otimes\G)^s$. If we introduce the notation $\tilde\C^r(\F):=\bigoplus_\Lambda i_{\Lambda *}(\F|_{V_\Lambda})$ for a sheaf $\F$ on $X\times X$, where the sum runs only over strictly increasing subsequences of $\Gamma$ of length $r+1$ (not subsequences of $\Gamma\sqcup\Gamma'$), then we have
\begin{align*}
&\G^{-3}=\tilde\C^0(\J/\J^2)\\
&\G^{-2}=\tilde\C^1(\J/\J^2)\oplus\tilde\C^0(\I_{\Delta_U}|_{X\times X})\oplus\tilde\C^0(\J/\J^2)\\
&\G^{-1}=\tilde\C^2(\J/\J^2)\oplus\tilde\C^1(\I_{\Delta_U}|_{X\times X})\oplus\tilde\C^1(\J/\J^2)\oplus\tilde\C^0(\mathcal O_{X\times X})\oplus\tilde\C^0(\I_{\Delta_U}/\I_{\Delta_U}^2|_{X\times X})\\
&\G^0=\tilde\C^3(\J/\J^2)\oplus\tilde\C^2(\I_{\Delta_U}|_{X\times X})\oplus\tilde\C^2(\J/\J^2)\oplus\tilde\C^1(\mathcal O_{X\times X})\oplus\tilde\C^1(\I_{\Delta_U}/\I_{\Delta_U}^2|_{X\times X})\oplus\tilde\C^0(\mathcal O_X).
\end{align*}
Thus $q_*(p^*\E\otimes\G)^s$ has the summand
\begin{align}\label{vekt2}
&q_*(p^*\E^{s}\otimes\tilde\C^3(\J/\J^2))\oplus q_*(p^*\E^{s}\otimes\tilde\C^2(\I_{\Delta_U}|_{X\times X}))\oplus q_*(p^*\E^{s}\otimes\tilde\C^2(\J/\J^2))\oplus\\&\nonumber
q_*(p^*\E^{s}\otimes\tilde\C^1(\mathcal O_{X\times X}))\oplus q_*(p^*\E^{s}\otimes\tilde\C^1(\I_{\Delta_U}/\I_{\Delta_U}^2|_{X\times X}))\oplus q_*(p^*\E^{s}\otimes\tilde\C^0(\mathcal O_X))\oplus\\&\nonumber
q_*(p^*\E^{s+1}\otimes
\tilde\C^2(\J/\J^2))\oplus q_*(p^*\E^{s+1}\otimes\tilde\C^1(\I_{\Delta_U}|_{X\times X}))\oplus q_*(p^*\E^{s+1}\otimes\tilde\C^1(\J/\J^2))\oplus\\&\nonumber
q_*(p^*\E^{s+1}\otimes\tilde\C^0(\mathcal O_{X\times X}))\oplus q_*(p^*\E^{s+1}\otimes\tilde\C^0(\I_{\Delta_U}/\I_{\Delta_U}^2|_{X\times X}))\oplus q_*(p^*\E^{s+2}\otimes
\tilde\C^1(\J/\J^2))\oplus\\&\nonumber
q_*(p^*\E^{s+2}\otimes\tilde\C^0(\I_{\Delta_U}|_{X\times X}))\oplus q_*(p^*\E^{s+2}\otimes\tilde\C^0(\J/\J^2))\oplus q_*(p^*\E^{s+3}\otimes
\tilde\C^0(\J/\J^2)).
\end{align}
The restrictions of $\J/\J^2$, $\I_{\Delta_U}|_{X\times X}$, $\mathcal O_{X\times X}$, $\I_{\Delta_U}/\I_{\Delta_U}^2|_{X\times X}$, and $\mathcal O_X$ to open sets of the form $V_\Lambda$ (with $\Lambda\subseteq\Gamma$) correspond to the $A_\Lambda/J_\Lambda\otimes_{k_\Lambda}A_\Lambda/J_\Lambda$-modules $J_\Lambda/J_\Lambda^2$, $I_\Lambda\otimes_{A_\Lambda\otimes_{k_\Lambda}A_\Lambda}(A_\Lambda/J_\Lambda\otimes_{k_\Lambda}A_\Lambda/J_\Lambda)$, $A_\Lambda/J_\Lambda\otimes_{k_\Lambda}A_\Lambda/J_\Lambda$, $I_\Lambda/I_\Lambda^2\otimes_{A_\Lambda\otimes_{k_\Lambda}A_\Lambda}(A_\Lambda/J_\Lambda\otimes_{k_\Lambda}A_\Lambda/J_\Lambda)$, and $A_\Lambda/J_\Lambda$, where $I_\Lambda\subseteq A_\Lambda\otimes_{k_\Lambda}A_\Lambda$ denotes the ideal of the diagonal. A morphism $\E^s\rightarrow q_*(p^*\E^t\otimes i_{\Lambda *}(\F|_{V_\Lambda}))$ with one of these sheaves $\F$ corresponds to a $m_i^t\times m_i^s$-matrix with entries in the modules mentioned, where $i$ is the minimal element in $\Lambda$, and where we use the restrictions of $\varphi_i^s$ bzw. $\varphi_i^t$ as trivializations of $\E^s$ and $\E^t$ on $X_\Lambda$.

The diagrams (\ref{univAti2}) and (\ref{vekt0}) show that the map $\E^s\rightarrow q_*(p^*\E\otimes\G)^s$ induced by the truncated Atiyah class factorizes over the summand (\ref{vekt2}), and that, in the above sense, it can be described by
\begin{equation*}
\left(0,0,0,0,0,(\mathbbm 1)_i,0,0,0,0,0,0,0,0,0\right),
\end{equation*}
where $\mathbbm 1$ denotes the identity matrix. Correspondingly, the morphism defined by (\ref{FormelAtiyah}) is given by
\begin{align*}
&\left(0,0,\left(
M^s_{ik}(\tilde M^s_{kj}\cdot\tilde M^s_{ji}-\tilde M^s_{ki})
\right)_{ijk},0,\left(
M^s_{ij}\cdot d\tilde M^s_{ji}
\right)_{ij},0,0,0,\right.\\&\left.\left(
-\xi M_{ij}^{s+1}(\tilde M_{ji}^{s+1}\cdot\tilde D_i^s-\tilde D_j^s\cdot\tilde M_{ji}^s)
\right)_{ij},0,\left(
-\xi d\tilde D_i^s
\right)_i,0,0,\left(
-\tilde D_i^{s+1}\cdot\tilde D_i^s
\right)_i,0\right).
\end{align*}
Once more, we write $\xi:=(-1)^s$.

We have to construct a homotopy between these two morphisms $\E\rightarrow q_*(p^*\E\otimes\G)$. This requires a map $h^s:\E^s\rightarrow q_*(p^*\E\otimes\G)^{s-1}$ in each degree. We define such a morphism, factorizing over the summand
\begin{align}\label{vekt3}
&
q_*(p^*\E^{s}\otimes
\tilde\C^2(\J/\J^2))\oplus q_*(p^*\E^{s}\otimes\tilde\C^1(\I_{\Delta_U}|_{X\times X}))\oplus
q_*(p^*\E^{s}\otimes\tilde\C^0(\mathcal O_{X\times X}))\oplus\\&\nonumber
q_*(p^*\E^{s+1}\otimes
\tilde\C^1(\J/\J^2))\oplus q_*(p^*\E^{s+1}\otimes\tilde\C^0(\I_{\Delta_U}|_{X\times X}))\oplus q_*(p^*\E^{s+2}\otimes
\tilde\C^0(\J/\J^2))
\end{align}
of $q_*(p^*\E\otimes\G)^{s-1}$, by
\begin{align*}
&\left(\left(
-\xi M^s_{ik}(\tilde M^s_{kj}\cdot\tilde M^s_{ji}-\tilde M^s_{ki})
\right)_{ijk},\left(
\xi((\tilde M_{ij}^s\otimes 1)(1\otimes\tilde M_{ji}^s)-(\tilde M_{ij}^s\cdot\tilde M_{ji}^s)\otimes 1)
\right)_{ij},\left(
-\xi\mathbbm 1
\right)_{i},\right.\\&\left.\left(
M_{ij}^{s+1}(\tilde M_{ji}^{s+1}\cdot\tilde D_i^s-\tilde D_j^s\cdot\tilde M_{ji}^s)
\right)_{ij},\left(
\tilde D_i^s\otimes 1-1\otimes\tilde D_i^s
\right)_i,\left(
\xi\tilde D_i^{s+1}\cdot\tilde D_i^s
\right)_i\right).
\end{align*}
Here, for a matrix $M$ with entries in $A_\Lambda$, we denote by $1\otimes M$ or $M\otimes 1$ the matrices with entries in $A_\Lambda\otimes_{k_\Lambda}A_\Lambda$ that result from $M$ by replacing each entry $m$ by $1\otimes m$ or $m\otimes 1$.

To conclude, we have to prove that in this way, we really get a homotopy between the morphisms mentioned. In the diagram
\begin{equation*}
\begin{xy}
\xymatrix{
&&q_*(p^*\E\otimes\G)^{s-1}\ar[dd]^{\hspace{1.5cm}{\mbox{\scriptsize$\left(\begin{matrix}\xi\check d&0&0&0&0&0\\-\xi\beta&\xi\check d&0&0&0&0\\\xi {\rm id}&0&0&0&0&0\\0&\xi\alpha&\xi\check d&0&0&0\\0&-\xi\pi&0&0&0&0\\0&0&-\xi\varepsilon&0&0&0\\
d_\E&0&0&-\xi\check d&0&0\\0&d_\E&0&-\xi\beta&-\xi\check d&0\\0&0&0&\xi {\rm id}&0&0\\0&0&d_\E&0&\xi\alpha&0\\0&0&0&0&-\xi\pi&0\\
0&0&0&d_\E&0&\xi\check d\\0&0&0&0&d_\E&-\xi\beta\\0&0&0&0&0&\xi {\rm id}\\
0&0&0&0&0&d_\E
\end{matrix}\right)$}}}\\\\
\E^s\ar[rruu]^{h^s}\ar[dd]^{d_\E}\ar@<0.1cm>[rr]^/-0.2cm/{At(\E)}\ar@<-0.1cm>[rr]_/-0.2cm/{(\ref{FormelAtiyah})}&&q_*(p^*\E\otimes\G)^{s}\\\\
\E^{s+1},\ar[rruu]_{h^{s+1}}
}\end{xy}
\end{equation*}
we only give the restriction of the differential of the complex $q_*(p^*\E\otimes\G)$ to the summand (\ref{vekt3}) of $q_*(p^*\E\otimes\G)^{s-1}$ which is relevant for our calculation, and which is mapped by the differential into the summand (\ref{vekt2}) of $q_*(p^*\E\otimes\G)^{s}$. Corresponding to the decomposition (\ref{vekt2}), we have to execute fifteen calculations. Fortunately, we have already finished the first, seventh, twelfth, and fifteenth one in the proof of Lemma \ref{l33} (as first, third, fifth, and seventh calculation there; the situations only differ by the factor $-\xi$). Furthermore, the calculations three, six, nine, and fourteen are trivial, and for the fifth and eleventh one, we only have to recall the identification $\Omega_{A_\Lambda|k_\Lambda}\cong I_\Lambda/I_\Lambda^2,da\mapsto 1\otimes a-a\otimes 1$.

Concerning the remaining summands of the decomposition (\ref{vekt2}), we proceed as in the proof of Lemma \ref{l33}. In particular, for each calculation, we fix a strictly increasing sequence $\Lambda$ of indices.

For the second summand, we find
\begin{align*}
&-\xi\beta(-\xi M^s_{ik}(\tilde M^s_{kj}\cdot\tilde M^s_{ji}-\tilde M^s_{ki}))+\xi\check d((\xi((\tilde M^s_{ij}\otimes 1)(1\otimes\tilde M^s_{ji})-(\tilde M^s_{ij}\cdot\tilde M^s_{ji})\otimes 1))_{ij})_{ijk}=\\&
-[-(\tilde M^s_{ik}\cdot\tilde M^s_{kj}\cdot\tilde M^s_{ji})\otimes 1+1\otimes(\tilde M^s_{ik}\cdot\tilde M^s_{kj}\cdot\tilde M^s_{ji})+(\tilde M^s_{ik}\cdot\tilde M^s_{ki})\otimes 1-1\otimes(\tilde M^s_{ik}\cdot\tilde M^s_{ki})]+\\&
[(M^s_{ij}\otimes 1)((\tilde M^s_{jk}\otimes 1)(1\otimes\tilde M^s_{kj})-(\tilde M^s_{jk}\cdot\tilde M^s_{kj})\otimes 1)(1\otimes M^s_{ji})-\\&(\tilde M^s_{ik}\otimes 1)(1\otimes\tilde M^s_{ki})+(\tilde M^s_{ik}\cdot\tilde M^s_{ki})\otimes 1+(\tilde M^s_{ij}\otimes 1)(1\otimes\tilde M^s_{ji})-(\tilde M^s_{ij}\cdot\tilde M^s_{ji})\otimes 1]=\\&
(\tilde M^s_{ik}\cdot\tilde M^s_{kj}\cdot\tilde M^s_{ji})\otimes 1-1\otimes(\tilde M^s_{ik}\cdot\tilde M^s_{kj}\cdot\tilde M^s_{ji})+1\otimes(\tilde M^s_{ik}\cdot\tilde M^s_{ki})+\\&
((\tilde M^s_{ij}\cdot\tilde M^s_{jk})\otimes 1)(1\otimes(\tilde M^s_{kj}\cdot\tilde M^s_{ji}))-((\tilde M^s_{ij}\cdot\tilde M^s_{jk}\cdot\tilde M^s_{kj})\otimes 1)(1\otimes \tilde M^s_{ji})-\\&
(\tilde M^s_{ik}\otimes 1)(1\otimes\tilde M^s_{ki})+(\tilde M^s_{ij}\otimes 1)(1\otimes\tilde M^s_{ji})-(\tilde M^s_{ij}\cdot\tilde M^s_{ji})\otimes 1=\\&
(\tilde M^s_{ik}\otimes 1-1\otimes\tilde M^s_{ik})(1\otimes(M^s_{kj}\cdot M^s_{ji}-M^s_{ki}))+\\&((M^s_{ij}\cdot M^s_{jk}\cdot M^s_{kj}-M^s_{ij})\otimes 1)(\tilde M^s_{ji}\otimes 1-1\otimes\tilde M^s_{ji})+\\&((M^s_{ij}\cdot M^s_{jk}-M^s_{ik})\otimes 1)(1\otimes(\tilde M^s_{kj}\cdot\tilde M^s_{ji})-(\tilde M^s_{kj}\cdot\tilde M^s_{ji})\otimes 1)=0.
\end{align*}

For the fourth summand, we get
\begin{align*}
&\xi\alpha(\xi((\tilde M^s_{ij}\otimes 1)(1\otimes\tilde M^s_{ji})-(\tilde M^s_{ij}\cdot\tilde M^s_{ji})\otimes 1))+\xi\check d((-\xi\mathbbm 1)_i)_{ij}=\\&
((M^s_{ij}\otimes 1)(1\otimes M^s_{ji})-(M^s_{ij}\cdot M^s_{ji})\otimes 1)+(-(M^s_{ij}\otimes 1)(1\otimes M^s_{ji})+\mathbbm 1)=0.
\end{align*}

For the eighth summand, we get
\begin{align*}
&d_\E(\xi((\tilde M_{ij}^s\otimes 1)(1\otimes\tilde M_{ji}^s)-(\tilde M_{ij}^s\cdot\tilde M_{ji}^s)\otimes 1))-\xi\beta(M_{ij}^{s+1}(\tilde M_{ji}^{s+1}\cdot\tilde D_i^s-\tilde D_j^s\cdot\tilde M_{ji}^s))-\\&\xi\check d((\tilde D_i^s\otimes 1-1\otimes\tilde D_i^s)_i)_{ij}+(-\xi((\tilde M_{ij}^{s+1}\otimes 1)(1\otimes\tilde M_{ji}^{s+1})-(\tilde M_{ij}^{s+1}\cdot\tilde M_{ji}^{s+1})\otimes 1))d_\E=\\&
\xi[(D_i^s\otimes 1)((\tilde M_{ij}^s\otimes 1)(1\otimes\tilde M_{ji}^s)-(\tilde M_{ij}^s\cdot\tilde M_{ji}^s)\otimes 1)-\\&
(\tilde M_{ij}^{s+1}\cdot\tilde M_{ji}^{s+1}\cdot\tilde D_i^s)\otimes 1+1\otimes(\tilde M_{ij}^{s+1}\cdot\tilde M_{ji}^{s+1}\cdot\tilde D_i^s)+(\tilde M_{ij}^{s+1}\cdot\tilde D_j^s\cdot\tilde M_{ji}^s)\otimes 1-\\&
1\otimes(\tilde M_{ij}^{s+1}\cdot\tilde D_j^s\cdot\tilde M_{ji}^s)-
(M_{ij}^{s+1}\otimes 1)(\tilde D_j^s\otimes 1-1\otimes\tilde D_j^s)(1\otimes M_{ji}^s)+\\&
\tilde D_i^s\otimes 1-1\otimes\tilde D_i^s-
((\tilde M_{ij}^{s+1}\otimes 1)(1\otimes\tilde M_{ji}^{s+1})-(\tilde M_{ij}^{s+1}\cdot\tilde M_{ji}^{s+1})\otimes 1)(1\otimes D^s_i)]=\\&
\xi[((\tilde D_i^s\cdot\tilde M_{ij}^s)\otimes 1)(1\otimes\tilde M_{ji}^s)-(\tilde D_i^s\cdot\tilde M_{ij}^{s}\cdot\tilde M_{ji}^{s})\otimes 1-(\tilde M_{ij}^{s+1}\cdot\tilde M_{ji}^{s+1}\cdot\tilde D_i^s)\otimes 1+\\&1\otimes(\tilde M_{ij}^{s+1}\cdot\tilde M_{ji}^{s+1}\cdot\tilde D_i^s)+(\tilde M_{ij}^{s+1}\cdot\tilde D_j^s\cdot\tilde M_{ji}^s)\otimes 1-1\otimes(\tilde M_{ij}^{s+1}\cdot\tilde D_j^s\cdot\tilde M_{ji}^s)-\\&
((\tilde M_{ij}^{s+1}\cdot\tilde D_j^s)\otimes 1)(1\otimes\tilde M_{ji}^s)+(\tilde M_{ij}^{s+1}\otimes 1)(1\otimes(\tilde D_j^s\cdot\tilde M_{ji}^s))+\tilde D_i^s\otimes 1-1\otimes\tilde D_i^s-\\&
(\tilde M_{ij}^{s+1}\otimes 1)(1\otimes(\tilde M_{ji}^{s+1}\cdot\tilde D^s_i))+((\tilde M_{ij}^{s+1}\cdot\tilde M_{ji}^{s+1})\otimes 1)(1\otimes\tilde D^s_i)]=\\&
\xi[((\mathbbm 1-M_{ij}^{s+1}\cdot M_{ji}^{s+1})\otimes 1)(\tilde D_i^s\otimes 1-1\otimes\tilde D_i^s)+\\&(\tilde M_{ij}^{s+1}\otimes1-1\otimes \tilde M_{ij}^{s+1})(1\otimes(D_j^s\cdot M_{ji}^s-M_{ji}^{s+1}\cdot D_i^s))+\\&((D_i^{s}\cdot M_{ij}^s-M_{ij}^{s+1}\cdot D_j^s)\otimes 1)(1\otimes\tilde M_{ji}^s-\tilde M_{ji}^s\otimes 1)]=0.
\end{align*}

For the tenth summand, we compute
\begin{align*}
&d_\E(-\xi\mathbbm 1)+\xi\alpha(\tilde D_i^s\otimes 1-1\otimes\tilde D_i^s)+(\xi\mathbbm 1)d_\E=\\&
\xi[-(D^s_i\otimes 1)\mathbbm 1+D_i^s\otimes 1-1\otimes D_i^s+\mathbbm 1(1\otimes D^s_i)]=0.
\end{align*}

Finally, we get
\begin{align*}
&d_\E(\tilde D_i^s\otimes 1-1\otimes\tilde D_i^s)-\xi\beta(\xi\tilde D_i^{s+1}\cdot\tilde D_i^s)+(\tilde D_i^{s+1}\otimes 1-1\otimes\tilde D_i^{s+1})d_\E=\\&
(D_i^{s+1}\otimes 1)(\tilde D_i^s\otimes 1-1\otimes\tilde D_i^s)-(\tilde D_i^{s+1}\cdot\tilde D_i^s)\otimes 1+1\otimes(\tilde D_i^{s+1}\cdot\tilde D_i^s)+\\&(\tilde D_i^{s+1}\otimes 1-1\otimes\tilde D_i^{s+1})(1\otimes D_i^s)=0
\end{align*}
for the thirteenth summand.
\hfill$\Box$
\end{Bew}

\subsection{The classical Atiyah class of a complex of vector bundles}

In the last subsection, we found a description of the truncated Atiyah class. There is a corresponding description of the classical class, found by \textit{Angéniol} and \textit{Lejeune-Jalabert} in \cite[Prop. II.1.5.2.1]{Ang}. Since the truncated class lifts the classical one, our formula allows us to recover the already known one.

As in the truncated case, morphisms $\E^s\rightarrow(\E\otimes\C(\Omega_X[1]))^s$ that factorize over the summand $\E^{s}\otimes\C^1(\Omega_X)\oplus\E^{s+1}\otimes\C^0(\Omega_X)$ of $(\E\otimes\C(\Omega_X[1]))^s$ can be constructed by giving matrices
\begin{equation}\label{FormelklasAti}
\left(\left(
M^s_{ij}\cdot d M^s_{ji}
\right)_{ij},\left(
(-1)^{s+1}d D_i^s
\right)_i\right).
\end{equation}

\begin{Kor}[Angéniol, Lejeune-Jalabert]
In situation \ref{sit}, the maps (\ref{FormelklasAti}) define a morphism of complexes $\E\rightarrow\E\otimes\C(\Omega_X[1])$. The induced morphism $\E\rightarrow\E\otimes\Omega_X[1]$ in the derived category is the classical Atiyah class of $\E$.
\end{Kor}

\begin{Bew}
The diagram
\begin{equation*}
\begin{xy}
\xymatrix{
&\E\otimes\C(\mathbb L_X[1])\ar[dd]&\E\otimes\mathbb L_X[1]\ar[dd]\ar[l]\\
\E\ar[ur]^/-0.2cm/{(\ref{FormelAtiyah})}\ar[dr]_/-0.2cm/{(\ref{FormelklasAti})}\\
&\E\otimes\C(\Omega_X[1])&\E\otimes\Omega_X[1]\ar[l]
}\end{xy}
\end{equation*}
is commutative by definition of the morphism $\mathbb L_X\rightarrow\Omega_X$ and the maps (\ref{FormelAtiyah}) and (\ref{FormelklasAti}).
Thus Lemma \ref{l23} gives the assertion of the corollary.
\hfill$\Box$
\end{Bew}

The formula given in \cite[Prop. II.1.5.2.1]{Ang} and our result (\ref{FormelklasAti}) differ by signs. Apart from the different conventions (that we discussed at the end of the introduction), there is another reason for this. Whereas we calculate with the complex $\E\otimes\C(\Omega_X[1])$ the complex $\C(\Omega_X\otimes\E)[1]$ is used in \cite{Ang}. There is a natural isomorphism between these complexes bringing further signs into play.

\section{The first truncated Chern class}

Via a trace map, the classical Atiyah class of a perfect complex gives rise to an element in the first cohomology group of the cotangent sheaf, called the first Chern class. Similarly, the truncated Atiyah class induces an element in the first cohomology group of truncated cotangent complex, the first truncated Chern class. This class will be investigated in this section. We will find out that the basic properties of the classical Chern class can also be proven for the truncated version. At the end of the section, we will discuss an application.

\subsection{The trace map}

First, we have to recall the trace map, introduced by \textit{Illusie} in \cite[Sect. I.8]{SGA6}. In \cite[Sect. V.3]{Ill}, he generalizes his construction in order to define Chern classes. We need the notion of a perfect complex. Recall that a complex $\E$ of sheaves on a scheme $X$ is said to be \emph{perfect} if every point of $X$ is contained in an open set $Y\subseteq X$ on which there exist a bounded complex $\F$ of finite free sheaves and a quasi-isomorphism $\F\rightarrow\E|_Y$ of complexes on $Y$.

\vspace{0.3cm}
Let $X$ be a scheme, let $\E$ be a perfect complex, and let $\G$ be a bounded complex of sheaves on $X$. In the derived category $\D({\rm Mod}(X))$, there are a natural isomorphism $R\H om(\E,\E\otimes^L \G)\xleftarrow{\sim}\E^\vee\otimes^L\E\otimes^L\G$ (see \cite[Cor. I.7.7]{SGA6}) and a map $\E^\vee\otimes^L\E\otimes^L\G\rightarrow\G$ which is induced by the contraction morphism $\E^\vee\otimes^L\E\rightarrow\mathcal O_X$ (see \cite[Lem. I.7.5]{SGA6}), where we use the abbreviation $\E^\vee=R\mathcal{H}om(\E,\mathcal O_X)$. The composition $R\H om(\E,\E\otimes^L \G)\rightarrow \G$ gives the trace map
\begin{equation*}
{\rm tr}:{\rm Hom}(\E,\E\otimes^L\G)\rightarrow{\rm Hom}(\mathcal O_X,\G)
\end{equation*}
after application of the functor ${\rm Hom}(\mathcal O_X,\_)$ and in view of the isomorphism ${\rm Hom}(\E,\E\otimes^L\G)\xrightarrow\sim{\rm Hom}(\mathcal O_X,R\H om(\E,\E\otimes^L \G))$ given by adjunction (see \cite[Lem. I.7.4]{SGA6}).

This definition of the trace map is complicated. The following lemma collects all its properties that are relevant for us.

\begin{Lemma}[Illusie]\label{l41}
Let $X$ be a scheme. Let $\E$ and $\F$ be perfect complexes, and let $\G$ and $\H$ be bounded complexes of sheaves on $X$. Let $\mu:\E\rightarrow\E\otimes^L\G$ be a morphism in $\D({\rm Mod}(X))$.
\begin{enumerate}[(i)]
\item{If $\E$ is a bounded complex of locally free coherent sheaves, and $\mu$ is a morphism of complexes, given by maps $\mu^{st}:\E^s\rightarrow\E^t\otimes\G^{s-t}$ of sheaves, then ${\rm tr}(\mu)$ is a morphism of complexes as well, given by $\sum_s (-1)^s {\rm tr}(\mu^{ss}):\mathcal O_X\rightarrow G^0$. The term ${\rm tr(\mu^{ss})}$ is to be understood in the following way: $\mu^{ss}$ corresponds to a map $\mathcal O_X\rightarrow(\E^s)^\vee\otimes\E^s\otimes\G^0$. This map has to be composed with the morphism $(\E^s)^\vee\otimes\E^s\otimes\G^0\rightarrow\G^0$ induced by the usual trace map $(\E^s)^\vee\otimes\E^s\rightarrow\mathcal O_X$.}
\item{If $\chi:\E\rightarrow\F$ is an isomorphism in $\D({\rm Mod}(X))$, and $\nu:=(\chi\otimes^L{\rm id})\circ\mu\circ\chi^{-1}:\F\rightarrow\F\otimes^L\G$ is the morphism induced by $\mu$, then ${\rm tr}(\mu)={\rm tr}(\nu)$.}
\item{If $\tau:\G\rightarrow\H$ is a morphism in $\D({\rm Mod}(X))$, then the diagram
\begin{equation*}
\begin{xy}
\xymatrix{
&&\G\ar[dd]^{\tau}\\
\mathcal O_X\ar[rru]^{{\rm tr}(\mu)}\ar[rrd]_{{\rm tr}((\E\otimes^L\tau)\circ\mu)}\\
&&\H
}\end{xy}
\end{equation*}
commutes.}
\end{enumerate}
\end{Lemma}

\begin{Bew}
Statement (i) is proved in \cite[I.8.1.2]{SGA6} in slightly weaker generality.
Statements (ii) and (iii) are discussed in \cite[V.3.8.11.1]{Ill} and \cite[V.3.7.4.1]{Ill}.
\hfill$\Box$
\end{Bew}

We are interested in the trace of the truncated Atiyah class.

\begin{Def}
Let $X$ be a separated $S$-scheme of finite type, and let $\E$ be a perfect complex on $X$. Then $c_{1,cl}(\E):={\rm tr}(At_{cl}(\E))\in H^1(X,\Omega_X)$ is said to be the \emph{first (classical) Chern class} of $\E$. If $X$ is flat with a smooth ambient space (over $S$), then $c_1(\E):={\rm tr}(At(\E))\in H^1(X,\mathbb L_X)$ is said to be the \emph{first truncated Chern class} of $\E$.
\end{Def}

From Lemma \ref{l23} and part (iii) of Lemma \ref{l41}, it follows that in the situation of the definition (with flat $X$ with smooth ambient space), the first classical Chern class is the composition $c_{1,cl}(\E):\mathcal O_X\xrightarrow{c_1(\E)}\mathbb L_X[1]\rightarrow\Omega_X[1]$.
Thus via the natural map $H^1(X,\mathbb L_X)\rightarrow H^1(X,\Omega_X)$, the truncated class lifts the classical one.

\subsection{The first truncated Chern class of a complex of vector bundles}

For bounded complexes of vector bundles, we can find an explicit formula for their first truncated Chern class. It can be derived easily from the results of the last section. Therefore, for this subsection, we return to the situation \ref{sit}.

In particular, we consider a flat $S$-scheme $X$ with a smooth ambient space $X\hookrightarrow U$ and a bounded complex $\E$ of locally free coherent sheaves, for which we have chosen a suitable cover $U=\bigcup_\Gamma U_i$ and trivializations of all components with transition maps $M^s_{ij}\in {\rm GL}(m_i^s,A_{ij}/J_{ij})$ and lifts $\tilde M^s_{ij}\in Mat(m^s_i\times m_i^s,A_{ij})$. For the local differentials $D^s_i\in Mat(m_i^{s+1}\times m_i^s,A_i/J_i)$, we have also chosen lifts $\tilde D^s_i\in Mat(m_i^{s+1}\times m_i^s,A_i)$.

\begin{Satz}\label{s42}
In situation \ref{sit}, the first truncated Chern class of $\E$ is given by the \v Cech cocycle
\begin{equation}\label{FormelChern}
\left(\left(
\sum_s (-1)^{s}{\rm tr}(M^s_{ik}(\tilde M^s_{kj}\cdot\tilde M^s_{ji}-\tilde M^s_{ki}))
\right)_{ijk},\left(
\sum_s (-1)^{s}{\rm tr}(M^s_{ij}\cdot d\tilde M^s_{ji})
\right)_{ij}\right)
\end{equation}
$\in H^0(X,\C^0(\mathbb L_X[1]))=H^0(X,\C^2(\J/\J^2))\oplus H^0(X,\C^1(\Omega_U|_X))$. More precisely, the formula (\ref{FormelChern}) defines a morphism $\mathcal O_X\rightarrow\C(\mathbb L_X[1])$, and the resulting diagram $\mathcal O_X\rightarrow\C(\mathbb L_X[1])\leftarrow\mathbb L_X[1]$ describes the first truncated Chern class of $\E$.
\end{Satz}

\begin{Bew}
The truncated Atiyah class of $\E$ is given by the diagram $\E\xrightarrow{\mu}\E\otimes\C(\mathbb L_X[1])\leftarrow\E\otimes\mathbb L_X[1]$, where $\mu$ is the map (\ref{FormelAtiyah}). Part (i) of Lemma \ref{l41} shows that the formula (\ref{FormelChern}) describes the trace of $\mu$. The assertion of the theorem follows from part (iii) of the same lemma.
\hfill$\Box$
\end{Bew}

As in the case of the Atiyah class, we can recover the already known formula for the classical Chern class from (\ref{FormelChern}). The statement follows from Theorem \ref{s42} because the truncated class lifts the classical one.

\begin{Kor}[Angéniol, Lejeune-Jalabert]
In situation \ref{sit}, the first classical Chern class of $\E$ is given by the \v Cech cocycle
\begin{equation*}
\left(
\sum_s (-1)^{s}{\rm tr}(M^s_{ij}\cdot dM^s_{ji})
\right)_{ij}
\end{equation*}
$\in H^0(X,\C^0(\Omega_X[1]))=H^0(X,\C^1(\Omega_X))$.
\end{Kor}

\subsection{Basic properties of first truncated Chern classes}

Many formulas are known that relate classical Chern classes of associated vector bundles. Here, we want to derive similar formulas for truncated Chern classes.

\begin{Satz}
Let $X$ be a flat $S$-scheme with a smooth ambient space. Then the first truncated Chern class gives rise to a group homomorphism $c_1:{\rm Pic}(X)\rightarrow\ H^1(X,\mathbb L_X)$.
\end{Satz}

\begin{Bew}
Let $\M$ and $\N$ be line bundles on $X$. By choosing a smooth ambient space, a suitable cover, and trivializations of the line bundles, we can achieve a situation as in the previous section, and we stick to the notation used there. We denote the transition maps of $\M$ and $\N$ by $M_{ij}$ and $N_{ij}$ and the chosen lifts by $\tilde M_{ij}$ and $\tilde N_{ij}$. With respect to the induced trivialization of $\M\otimes\N$, the transition maps are given by $M_{ij}\cdot N_{ij}$, and we can choose the lifts $\tilde M_{ij}\cdot\tilde N_{ij}$.

By (\ref{FormelChern}), the first truncated Chern classes of $\M$ and $\N$ are given by
\begin{equation*}
\left(\left(M_{ik}(\tilde M_{kj}\cdot\tilde M_{ji}-\tilde M_{ki})\right)_{ijk},\left(M_{ij}\cdot d\tilde M_{ji}\right)_{ij}\right)
\end{equation*}
and
\begin{equation*}
\left(\left(N_{ik}(\tilde N_{kj}\cdot\tilde N_{ji}-\tilde N_{ki})\right)_{ijk},\left(N_{ij}\cdot d\tilde N_{ji}\right)_{ij}\right)
\end{equation*}
and those of $\M\otimes\N$ by
\begin{equation*}
\left(\left(M_{ik}\cdot N_{ik}(\tilde M_{kj}\cdot\tilde N_{kj}\cdot\tilde M_{ji}\cdot\tilde N_{ji}-\tilde M_{ki}\cdot\tilde N_{ki})\right)_{ijk},\left(M_{ij}\cdot N_{ij}\cdot d(\tilde M_{ji}\cdot\tilde N_{ji})\right)_{ij}\right).
\end{equation*}

To prove the theorem, we will show that the latter cocycle is the sum of the other two. We have
\begin{align*}
&M_{ik}(\tilde M_{kj}\cdot\tilde M_{ji}-\tilde M_{ki})+N_{ik}(\tilde N_{kj}\cdot\tilde N_{ji}-\tilde N_{ki})-\\&
M_{ik}\cdot N_{ik}(\tilde M_{kj}\cdot\tilde N_{kj}\cdot\tilde M_{ji}\cdot\tilde N_{ji}-\tilde M_{ki}\cdot\tilde N_{ki})=\\&
(\tilde M_{ik}\cdot\tilde M_{ki}-1)(\tilde N_{ik}\cdot\tilde N_{ki}-1)-\\&
(\tilde M_{ik}\cdot\tilde M_{kj}\cdot\tilde M_{ji}-1)(\tilde N_{ik}\cdot\tilde N_{kj}\cdot\tilde N_{ji}-1)=0
\end{align*}
(because the product of two elements of $J_{ijk}$ lies in $J_{ijk}^2$) and
\begin{align*}
&M_{ij}\cdot d\tilde M_{ji}+N_{ij}\cdot d\tilde N_{ji}=\\&M_{ij}\cdot N_{ij}\cdot N_{ji}\cdot d\tilde M_{ji}+M_{ij}\cdot N_{ij}\cdot M_{ji}\cdot d\tilde N_{ji}=\\&M_{ij}\cdot N_{ij}\cdot d(\tilde M_{ji}\cdot\tilde N_{ji}).
\end{align*}
\hfill$\Box$
\end{Bew}

We can associate to every bounded complex $\E$ of locally free coherent sheaves on a scheme $X$ a line bundle, called its determinant ${\rm det}(\E):=\bigotimes_s {\rm det}(\E^s)^{(-1)^s}$. The following main theorem of this section asserts that the first truncated Chern class of $\E$ only depends on its determinant.

\begin{Satz}\label{S45}
Let $X$ be a flat $S$-scheme with a smooth ambient space. Let $\E$ be a bounded complex of locally free coherent sheaves on $X$. Then for the first truncated Chern class, we have
\begin{equation*}
c_1(\E)=c_1({\rm det}(\E)).
\end{equation*}
\end{Satz}

\begin{Bew}
By choosing a smooth ambient space, a suitable cover, and trivializations, we can assume that we are in situation \ref{sit}. The formula (\ref{FormelChern}) shows that $c_1(\E)=\sum_s(-1)^s c_1(\E^s)$, and because of the additivity of the first truncated Chern class on the Picard group, we have $c_1({\rm det}(\E))=c_1(\bigotimes_s {\rm det}(\E^s)^{(-1)^s})=\sum_s(-1)^s c_1({\rm det}(\E^s))$. Thus we may assume that $\E$ is a single locally free coherent sheaf.

For the transition maps $M_{ij}$, we have chosen lifts $\tilde M_{ij}$. For the determinant of $\E$, there exist trivializations, with respect to which the transition maps are ${\rm det}(M_{ij})$. We choose the lifts ${\rm det}(\tilde M_{ij})$.

By (\ref{FormelChern}), the first Chern classes of $\E$ and ${\rm det}(\E)$ are given by 
\begin{equation*}
\left(\left({\rm tr}(M_{ik}(\tilde M_{kj}\cdot\tilde M_{ji}-\tilde M_{ki}))\right)_{ijk},\left(\sum_{u,v}(M_{ij})_{uv}\cdot d(\tilde M_{ji})_{vu}\right)_{ij}\right)
\end{equation*}
and
\begin{equation*}
\left(\left({\rm det}(M_{ik})({\rm det}(\tilde M_{kj})\cdot {\rm det}(\tilde M_{ji})-{\rm det}(\tilde M_{ki}))\right)_{ijk},\left({\rm det}(M_{ij})\cdot d({\rm det}(\tilde M_{ji}))\right)_{ij}\right).
\end{equation*}

We will prove that these cocycles coincide. To this end, we need the following simple result. Let $A$ be a ring with ideal $J$ and $M\in {\rm Mat}(m\times m,A)$ a matrix for which the entries of $M-\mathbbm 1$ all lie in $J$. Then in $A/J^2$, we have
\begin{align*}
&{\rm det}(M)=\sum_{\sigma\in S_m}sgn(\sigma)\prod_u M_{u,\sigma(u)}=\prod_u M_{uu}=\prod_u (1+(M_{uu}-1))=\\&1+\sum_u(M_{uu}-1)={\rm tr}(M)-(m-1).
\end{align*}

We use this result to prove the first of the required calculations:
\begin{align*}
&{\rm tr}(M_{ik}(\tilde M_{kj}\cdot\tilde M_{ji}-\tilde M_{ki}))={\rm tr}(\tilde M_{ik}\cdot\tilde M_{kj}\cdot\tilde M_{ji})-{\rm tr}(\tilde M_{ik}\cdot\tilde M_{ki})=\\&{\rm det}(\tilde M_{ik}\cdot\tilde M_{kj}\cdot\tilde M_{ji})-{\rm det}(\tilde M_{ik}\cdot\tilde M_{ki})={\rm det}(M_{ik})({\rm det}(\tilde M_{kj})\cdot {\rm det}(\tilde M_{ji})-{\rm det}(\tilde M_{ki})).
\end{align*}

For an algebra $A$ over a ring $k$ and a matrix $M\in {\rm Mat}(m\times m,A)$, we denote by $M^{uv}$ the matrix that results from $M$ by replacing the entry $(u,v)$ by $1$ and the other entries in row $u$ or column $v$ by $0$. It suffices to prove the formula
\begin{equation}\label{vekt1}
\sum_{u,v}{\rm det}(M^{vu})\cdot d M_{vu}=d({\rm det}(M))
\end{equation}
in $\Omega_{A|k}$ in this situation, for in our case $M=\tilde M_{ji}$, this formula implies the required equality when we use Cramer's rule:
\begin{equation*}
\sum_{u,v}(M_{ij})_{uv}\cdot d(\tilde M_{ji})_{vu}=\sum_{u,v}{\rm det}(M_{ij})\cdot {\rm det}((M_{ji})^{vu})\cdot d(\tilde M_{ji})_{vu}={\rm det}(M_{ij})\cdot d({\rm det}(\tilde M_{ji})).
\end{equation*}

We prove (\ref{vekt1}) by induction with respect to $m$, where the statement is clear for $m=1$. For $m>1$, we have
\begin{align*}
&d({\rm det}(M))=d(\sum_n M_{1n}\cdot {\rm det}(M^{1n}))=\\&\sum_n {\rm det}(M^{1n})\cdot d(M_{1n})+\sum_n M_{1n}\cdot d({\rm det}(M^{1n}))=\\&
\sum_n {\rm det}(M^{1n})\cdot d(M_{1n})+\sum_n M_{1n}\sum_{u\neq n,v\neq1}{\rm det}((M^{1n})^{vu})\cdot d(M_{vu})=\\&
\sum_u {\rm det}(M^{1u})\cdot d(M_{1u})+\sum_{v\neq1,u}(\sum_{n\neq u}M_{1n}\cdot {\rm det}((M^{vu})^{1n}))d(M_{vu})=\\&
\sum_{u,v}{\rm det}(M^{vu})\cdot d M_{vu},
\end{align*}
where we use Laplace's formula for the first and the last equality. To obtain the third equality, the induction hypothesis is applied to matrices that result from $M$ by deleting one row and one column.
\hfill$\Box$
\end{Bew}

To explain the following theorem, we have to discuss the rank of a bounded complex $\E$ of finite locally free sheaves on a scheme $X$. For a point $x\in X$, the number ${\rm rk}_x(\E):=\sum_s (-1)^s {\rm rk}_x \E^s$ is called the rank of $\E$ in $x$. If this number is the same for all points of $X$, we say that $\E$ is of constant rank. If $X$ is connected, then $\E$ obviously is of constant rank.

\begin{Satz}
Let $X$ be a flat $S$-scheme with smooth ambient space. Let $\E$, $\F$, and $\G$ be bounded complexes of locally free coherent sheaves on $X$. Then for the first truncated Chern classes, we have
\begin{equation}\label{Chern1}
c_1(\E^\vee)=-c_1(\E)
\end{equation}
\begin{equation}\label{Chern2}
c_1(\E\oplus\F)=c_1(\E)+c_1(\F).
\end{equation}
If there exists a short exact sequence $0\rightarrow\E\rightarrow\G\rightarrow\F\rightarrow 0$, then
\begin{equation}\label{Chern3}
c_1(\G)=c_1(\E)+c_1(\F).
\end{equation}
If $\E$ and $\F$ are of constant rank, then
\begin{equation}\label{Chern4}
c_1(\E\otimes\F)={\rm rk}(\F)\cdot c_1(\E)+{\rm rk}(\E)\cdot c_1(\F)
\end{equation}
\begin{equation}\label{Chern5}
c_1(\mathcal{H}om(\E,\F))=-{\rm rk}(\F)\cdot c_1(\E)+{\rm rk}(\E)\cdot c_1(\F).
\end{equation}
\end{Satz}

\begin{Bew}
To prove (\ref{Chern1}), we use the previous two theorems and find $c_1(\E^\vee)=c_1({\rm det}(\E^\vee))=c_1({\rm det}(\E)^\vee)=-c_1({\rm det}(\E))=-c_1(\E)$.

The equation (\ref{Chern3}) follows analogously from the calculation $c_1(\G)=c_1({\rm det}(\G))=c_1({\rm det}(\E)\otimes {\rm det}(\F))=c_1({\rm det}(\E))+c_1({\rm det}(\F))=c_1(\E)+c_1(\F)$, and it implies (\ref{Chern2}) as a special case.

Again using the previous theorems, we see that it suffices to prove
\begin{equation}\label{perf1}
{\rm det}(\E\otimes\F)={\rm det}(\E)^{\otimes {\rm rk}(\F)}\otimes {\rm det}(\F)^{\otimes {\rm rk}(\E)}
\end{equation}
to obtain (\ref{Chern4}). Since it is more difficult to prove this formula than its analogues above, we recall its proof. To this end, we can restrict to the connected components of $X$. Thus we can assume that the components of $\E$ and $\F$ are globally of constant rank. To prove (\ref{perf1}) for vector bundles $\E$ and $\F$, we fix trivializations of these sheaves on a suitable cover, and we compare the induced transition maps of the line bundles on both sides of the equation. For arbitrary square matrices $M$ and $N$ of size $m$ and $n$ with entries in a ring, however, we really have ${\rm det}(M\otimes N)={\rm det}((M\otimes\mathbbm 1)(\mathbbm 1\otimes N))={\rm det}(M\otimes\mathbbm 1)\cdot {\rm det}(\mathbbm 1\otimes N)={\rm det}(M)^n\cdot {\rm det}(N)^m$.

From this, we can deduce the statement in the general case of complexes:
\begin{align*}
&{\rm det}(\E\otimes\F)=\bigotimes_u {\rm det}((\E\otimes\F)^u)^{(-1)^u}=\bigotimes_u {\rm det}(\bigoplus_{s+t=u} \E^s\otimes\F^t)^{(-1)^u}=\\&
\bigotimes_{s,t}{\rm det}(\E^s\otimes\F^{t})^{(-1)^{s+t}}=\bigotimes_{s,t}({\rm det}(\E^s)^{\otimes {\rm rk}(\F^t)}\otimes {\rm det}(\F^t)^{\otimes {\rm rk}(\E^s)})^{(-1)^{s+t}}=\\&
\bigotimes_s {\rm det}(\E^s)^{\otimes(\sum_t(-1)^{s+t}{\rm rk}(\F^t))}\otimes\bigotimes_t {\rm det}(\F^t)^{\otimes(\sum_s(-1)^{s+t}{\rm rk}(\E^s))}=\\&
\bigotimes_s {\rm det}(\E^s)^{\otimes(-1)^{s}{\rm rk}(\F)}\otimes\bigotimes_t {\rm det}(\F^t)^{\otimes(-1)^{t}{\rm rk}(\E)}=
{\rm det}(\E)^{\otimes {\rm rk}(\F)}\otimes {\rm det}(\F)^{\otimes {\rm rk}(\E)}.
\end{align*}

Finally, the statement (\ref{Chern5}) follows from (\ref{Chern1}) and (\ref{Chern4}) if we use the isomorphism $\mathcal{H}om(\E,\F)\cong\E^\vee\otimes\F$.
\hfill$\Box$
\end{Bew}

\subsection{The first truncated Chern class of a perfect complex}

We defined the first truncated Chern class for all perfect complexes, but only considered it in the case of bounded complexes of vector bundles up to now. The reason is that only for this smaller class of complexes, we have the formula (\ref{FormelChern}) at our disposal. Although perfect complexes are locally isomorphic to complexes of this form, we cannot control their (globally defined) first truncated Chern classes with our methods. In order to derive statements as in the previous subsection for arbitrary perfect complexes, technically more complicated methods as in \cite[Chapt. V]{Ill} would probably be necessary.

Nevertheless, many schemes have the property that every perfect complex is even globally isomorphic to a bounded complex of vector bundles. In this case, our results can immediately be extended to perfect complexes. We will formulate these slightly more general statements in this subsection.

In this context, we will speak about determinants of perfect complexes, see \cite{KM}. We will consider the determinant of a perfect complex only as an element in the Picard group. The following lemma collects the most important properties (see \cite[Def. 4, Thm. 2, Sect. before Prop. 6]{KM}).

\begin{Lemma}[Knudsen, Mumford]
Let $X$ be a scheme, $\E$, $\F$, and $\G$ be perfect complexes on $X$. If $\E$ is a bounded complex of locally free coherent sheaves, then ${\rm det}(\E)=\bigotimes_s {\rm det}(\E^s)^{(-1)^s}$. If $\E$ and $\F$ are isomorphic in $\D({\rm Mod}(X))$, then ${\rm det}(\E)={\rm det}(\F)$. If there exists a distinguished triangle $\E\rightarrow\F\rightarrow\G\rightarrow\E[1]$, then ${\rm det}(\F)={\rm det}(\E)\otimes {\rm det}(\G)$.
\end{Lemma}

For bounded complexes of finite locally free sheaves, we introduced the rank in a point. Perfect complexes are locally isomorphic to such complexes. Thus the rank in a point is also defined for them. (To see that this notion is well-defined, we have to understand that two bounded complexes $\E$ and $\F$ of free finite modules which are isomorphic in the derived category of a local ring (the local ring in the point considered) have the same rank. After reduction to the case of a quasi-isomorphism $\E\rightarrow\F$, we consider its mapping cone $\G$. By construction of the cone, it suffices to prove ${\rm rk}(\G)=0$. Since $\G$ is an exact complex of free finite modules, the images of all differentials are projective and thus free (the ring is local). Hence $\G$ is split and exact.) All perfect complexes on connected schemes are of constant rank.

Now, we can easily generalize the results of the previous subsection (use (\ref{Cherndp}) for the proof of (\ref{Chern3p})).

\begin{Kor}\label{s48}
Let $X$ be a flat $S$-scheme with smooth ambient space. Let $\E$, $\F$, and $\G$ be perfect complexes on $X$ which are isomorphic to bounded complexes of locally free coherent sheaves in $\D({\rm Mod}(X))$. Then for the first truncated Chern classes, we have
\begin{equation}\label{Cherndp}
c_1(\E)=c_1({\rm det}(\E))
\end{equation}
\begin{equation*}
c_1(\E^\vee)=-c_1(\E)
\end{equation*}
\begin{equation*}
c_1(\E\oplus\F)=c_1(\E)+c_1(\F).
\end{equation*}
If there exists a distinguished triangle $\E\rightarrow\F\rightarrow\G\rightarrow\E[1]$, then
\begin{equation}\label{Chern3p}
c_1(\F)=c_1(\E)+c_1(\G).
\end{equation}
If $\E$ and $\F$ are of constant rank, then
\begin{equation*}
c_1(\E\otimes^L\F)={\rm rk}(\F)\cdot c_1(\E)+{\rm rk}(\E)\cdot c_1(\F)
\end{equation*}
\begin{equation*}
c_1(R\mathcal{H}om(\E,\F))=-{\rm rk}(\F)\cdot c_1(\E)+{\rm rk}(\E)\cdot c_1(\F).
\end{equation*}
\end{Kor}

If all perfect complexes on our scheme $X$ are isomorphic to bounded complexes of vector bundles, then the assumption of Corollary \ref{s48} is trivially satisfied. Therefore, we are interested in situations where this is the case. We will cite an interesting theorem due to \textit{Illusie} in this context. It states that the large class of divisorial schemes has the property mentioned. Recall that a separated quasi-compact scheme $X$ is said to be \emph{divisorial} if the open sets $X_s:=\{x\in X: s(x)\neq 0\}$ for arbitrary global sections $s$ of arbitrary line bundles on $X$ form a basis of the topology.

\vspace{0.3cm}
Let $A$ be a noetherian ring and $X$ a quasi-projective scheme over ${\rm Spec}(A)$. We show that in this situation, $X$ is divisorial. Let $x\in X$, and let $Y\subseteq X$ be an open neighbourhood of $x$. If we denote the ideal sheaf of $X\backslash Y$ (with the induced reduced scheme structure) by $\I$, and if $\M$ is an ample line bundle, then there exists a natural number $N$ for which $\I\otimes\M^{\otimes N}$ is globally generated. In particular, there is a global section $s$ with $s(x)\neq 0$. Via the inclusion $\I\subseteq\mathcal O_X$, we can consider it als a global section of the line bundle $\M^{\otimes N}$. This section has the required property $x\in X_s\subseteq Y$.

\begin{Satz}[Illusie]
Let $X$ be a separated, quasi-compact, and divisorial scheme (e.g. a quasi-projective scheme over a noetherian ring). Then every perfect complex on $X$ is isomorphic to a bounded complex of locally free coherent sheaves in $\D({\rm Mod}(X))$.
\end{Satz}

\begin{Bew}
A proof can be found in \cite[Prop. 2.3.1]{TT}. Originally, the theorem was proven in \cite[Chapt. II]{SGA6}.
\hfill$\Box$
\end{Bew}

If $S$ is quasi-projective over an affine noetherian scheme, and $X$ is flat and quasi-projective over $S$, then the theorem shows that $X$ is divisorial, and our Corollary \ref{s48} holds true for all perfect complexes. In this situation, the commutative diagram
\begin{equation*}
\begin{xy}
\xymatrix{
K({\rm Perf}(X))\ar[rr]^{\rm det}\ar[rd]_{c_1}&&{\rm Pic}(X)\ar[dl]^{c_1}\\
&H^1(X,\mathbb L_X)
}\end{xy}
\end{equation*}
of abelian groups summarizes the main results of this section.
Here, ${\rm Perf}(X)$ denotes the triangulated (see \cite[Prop. I.4.10]{SGA6}) category of perfect complexes on $X$, and $K({\rm Perf}(X))$ is its Grothendieck group.

\subsection{An application}

In this final subsection, we want to discuss a consequence of our Theorem \ref{S45} concerning the paper \cite{HT}.

We consider a closed immersion $i:X\hookrightarrow Y$ of flat $S$-schemes with smooth ambient spaces, and we assume that the ideal sheaf $\I$ of this immersion satisfies $\I^2=0$. We are interested in the question whether given perfect complexes $\E$ on $X$ extend to perfect complexes $\F$ on $Y$ in the sense that $Li^*(\F)\cong\E$. We call such extensions deformations of $\E$.

It turns out that the truncated Atiyah class of $\E$ can be used to give an answer to this question. Furthermore, the truncated Kodaira--Spencer class $\kappa(i):\mathbb L_{X}\rightarrow \I[1]$ (see \cite[Def. 2.7]{HT}) is useful --- note that $\I$ can be considered as sheaf on $X$ because of $\I^2=0$.

\begin{Satz}[Huybrechts, Thomas]
Let $X$ and $Y$ be flat $S$-schemes with smooth ambient spaces and $i:X\hookrightarrow Y$ a closed immersion whose ideal sheaf $\I$ satisfies $\I^2=0$. Let $\E$ be a perfect complex on $X$. Then the following statements are equivalent:
\begin{enumerate}[(i)]
\item{There exists a perfect complex $\F$ on $Y$ with $Li^*(\F)\cong\E$.}
\item{The map $\overline{\omega}(\E):=(\E\otimes^L\kappa(i)[1])\circ At(\E):\E\rightarrow\E\otimes^L\I[2]$ is the zero morphism.}
\end{enumerate}
\end{Satz}

\begin{Bew}
The theorem is proven in \cite[Cor. 3.4]{HT}.
\hfill$\Box$
\end{Bew}

For a perfect complex $\E$ on $X$, the morphism $\overline{\omega}(\E)$ is called the obstruction class of $\E$. We can use our results to relate the obstruction class of a perfect class to the obstruction class of its determinant. In this way, we give a more natural proof for a statement that \textit{Huybrechts} and \textit{Thomas} use in \cite{HT}.

\begin{Satz}
Let $X$ and $Y$ be flat $S$-schemes with smooth ambient spaces and $i:X\hookrightarrow Y$ a closed immersion whose ideal sheaf $\I$ satisfies $\I^2=0$. Let $\E$ be a perfect complex on $X$ that is isomorphic to a bounded complex of locally free coherent sheaves in $\D({\rm Mod}(X))$ (e.g. an arbitrary perfect complex if $X$ is divisorial). Then
\begin{equation*}
\overline{\omega}({\rm det}(\E))={\rm tr}(\overline{\omega}(\E)).
\end{equation*}
\end{Satz}

\begin{Bew}
For the equality, we make use of the natural isomorphism ${\rm Hom}({\rm det}(\E),{\rm det}(\E)\otimes^L\I[2])\cong{\rm Hom}(\mathcal O_X,\I[2])$. The precise formulation of the formula of the theorem is
\begin{equation*}
{\rm tr}(\overline{\omega}({\rm det}(\E)))={\rm tr}(\overline{\omega}(\E)).
\end{equation*}

We have ${\rm tr}(\overline{\omega}(\E))={\rm tr}((\E\otimes^L\kappa(i)[1])\circ At(\E))=\kappa(i)[1]\circ{\rm tr}(At(\E))=\kappa(i)[1]\circ c_1(\E)$, where we use part (iii) of Lemma \ref{l41} for the second equality. Correspondingly, we have ${\rm tr}(\overline{\omega}({\rm det}(\E)))=\kappa(i)[1]\circ c_1({\rm det}(\E))$. An application of formula (\ref{Cherndp}) finishes the proof.
\hfill$\Box$
\end{Bew}

\end{document}